\documentclass[a4paper,11pt]{article}
\usepackage{amsmath, amsfonts, amscd, amssymb, amsthm, enumerate, xypic}

\def\bfB{\mathbf{B}}

\DeclareMathOperator{\Kbar}{\overline{\mathbb{K}}}
\DeclareMathOperator{\id}{\operatorname{id}}
\DeclareMathOperator{\simt}{\underset{3}{\sim}}
\DeclareMathOperator{\com}{\operatorname{com}}
\DeclareMathOperator{\card}{\#\,}
\DeclareMathOperator{\Mat}{\operatorname{M}}
\DeclareMathOperator{\End}{\operatorname{End}}
\DeclareMathOperator{\Diag}{\operatorname{Diag}}
\DeclareMathOperator{\NT}{\operatorname{NT}}
\DeclareMathOperator{\GL}{\operatorname{GL}}
\DeclareMathOperator{\Ker}{\operatorname{Ker}}
\DeclareMathOperator{\Vect}{\operatorname{span}}
\DeclareMathOperator{\im}{\operatorname{Im}}
\DeclareMathOperator{\tr}{\operatorname{tr}}
\DeclareMathOperator{\car}{\operatorname{char}}
\DeclareMathOperator{\Sp}{\operatorname{Sp}}
\DeclareMathOperator{\rk}{\operatorname{rk}}
\DeclareMathOperator{\codim}{\operatorname{codim}}
\renewcommand{\setminus}{\smallsetminus}

\def\K{\mathbb{K}}

\def\C{\mathbb{C}}

\def\calF{\mathcal{F}}
\def\calG{\mathcal{G}}
\def\calH{\mathcal{H}}
\def\calI{\mathcal{I}}

\def\calT{\mathcal{T}}

\def\lcro{\mathopen{[\![}}
\def\rcro{\mathclose{]\!]}}

\theoremstyle{definition}
\newtheorem{Def}{Definition}
\newtheorem{Not}[Def]{Notation}

\theoremstyle{plain}
\newtheorem{theo}{Theorem}
\newtheorem{prop}[theo]{Proposition}
\newtheorem{cor}[theo]{Corollary}
\newtheorem{lemme}[theo]{Lemma}
\newtheorem{claim}{Claim}

\theoremstyle{plain}
\newtheorem{conj}{Conjecture}

\theoremstyle{remark}
\newtheorem{Rems}{Remarks}
\newtheorem{Rem}[Rems]{Remark}
\newtheorem{ex}[Rems]{Example}

\title{Spaces of matrices with a sole eigenvalue}
\author{Cl\'ement de Seguins Pazzis\footnote{Lyc\'ee Priv\'e Sainte-Genevi\`eve, 2, rue
de l'\'Ecole des Postes, 78029 Versailles Cedex, FRANCE.}
\footnote{e-mail address: dsp.prof@gmail.com}}

\begin{document}

\thispagestyle{plain}

\maketitle
\begin{abstract}
Let $\K$ be an arbitrary (commutative) field and $\overline{\K}$ be an algebraic closure of it.
Let $V$ be a linear subspace of $\Mat_n(\K)$, with $n \geq 3$. We show that if
every matrix of $V$ has at most one eigenvalue in $\K$, then $\dim V \leq 1+\binom{n}{2}$. \\
If every matrix of $V$ has a sole eigenvalue in $\overline{\K}$ and $\dim V=1+\binom{n}{2}$,
we show that $V$ is similar to the space of all upper-triangular matrices with equal diagonal entries,
except if $n=3$ and $\K$ has characteristic $3$, or if $n=4$ and $\K$ has characteristic $2$.
In both of those special cases, we classify the exceptional solutions up to similarity.
\end{abstract}



\section{Introduction}\label{pureintro}

In this article, we let $\K$ be an arbitrary (commutative) field,
and we choose an algebraic closure $\Kbar$ of it.
 We denote by
$\Mat_n(\K)$ the algebra of square matrices with $n$ rows
and entries in $\K$, and by $\GL_n(\K)$ its group of invertible elements.
We also denote by $\Mat_{n,p}(\K)$ the vector space of matrices with $n$ rows, $p$ columns and entries in $\K$.
The transpose of a matrix $M$ is denoted by $M^T$. \\
Two linear subspaces $V$ and $W$ of $\Mat_n(\K)$ are called \textbf{similar}, and we write $V \sim W$,
if $W=PVP^{-1}$ for some $P \in \GL_n(\K)$ (i.e., $V$ and $W$ represent, in a change of basis, the same set of endomorphisms of an
$n$-dimensional vector space).
For $M \in \Mat_n(\K)$, we denote by $\Sp_\K(M)$ its spectrum in $\K$, i.e., its set of eigenvalues in the field $\K$,
and by $\tr(M)$ its trace. \\
For $(a_1,\dots,a_n)\in \K^n$, we denote by $\Diag(a_1,\dots,a_n)$ the diagonal matrix of $\Mat_n(\K)$ with diagonal entries $a_1,\dots,a_n$. \\
Given two integers $a$ and $b$, we write $a \mid b$ if $a$ divides $b$, and
$a \wedge b=1$ if $a$ is prime with $b$.

\paragraph{}
Linear spaces of square matrices with conditions on their spectrum have been
the topic of quite a few papers in the past decades. The first important results
can be traced back to Gerstenhaber, who proved the following result in the case $\card \K \geq n$
(see also \cite{Mathes} for a simplified proof and an extension to the case $\card \K \geq 3$, and
\cite{Serezhkin} for a proof with no restriction on the field):

\begin{theo}[Gerstenhaber, Serezhkin]\label{Gerstentheo}
Let $V$ be a linear subspace of $\Mat_n(\K)$ in which every matrix is nilpotent. \\
Then $\dim V \leq \binom{n}{2}$. If equality holds, then $V$ is similar to the subspace
$\NT_n(\K)$ of strictly upper-triangular matrices.
\end{theo}

\paragraph{}
In \cite{Omladic},  Omladi\v c and \v Semrl consider the following general problem:
given $k \in \lcro 1,n-1\rcro$, determine the maximal dimension for a linear subspace of $\Mat_n(\K)$
in which every matrix has at most $k$ eigenvalues in $\Kbar$, and classify the subspaces with the maximal dimension.
They solved the problem in the special case $\K=\C$, for $k=1$, $k=2$ and $n$ odd, and $k=n-1$ with the additional
condition that the subspaces under consideration contain a matrix which has exactly $n-1$ distinct eigenvalues in $\Kbar$.
In the subsequent \cite{Loewy}, Loewy and Radwan considered mainly the ``upper bound" component of the problem,
and extended Omladi\v c and \v Semrl's results to an arbitrary field of characteristic $0$, whilst
solving the additional cases $k=3$ and $k=n-1$.

\paragraph{}
In this paper, we tackle the case $k=1$ for an arbitrary field
and extend the upper bound on the dimension to a larger class of subspaces.
Let us start with a few definitions.

\begin{Def}
Let $V$ be a linear subspace of $\Mat_n(\K)$. \\
We say that $V$ is a \textbf{$1$-spec subspace} when $\card \Sp_\K(M) \leq 1$ for every $M \in V$. \\
We say that $V$ is a \textbf{$\overline{1}$-spec subspace} when $\card \Sp_{\Kbar}(M)=1$ for every $M \in V$. \\
We say that $V$ has a \textbf{trivial spectrum} when $\Sp_\K(M) \subset \{0\}$ for every $M \in V$. \\
We say that $V$ is \textbf{nilpotent} when all its elements are nilpotent matrices, i.e., $\Sp_{\Kbar}(M)=\{0\}$ for every $M \in V$.
\end{Def}

Spaces with a trivial spectrum are linked to the affine subspaces of non-singular matrices of $\Mat_n(\K)$.
In \cite{Quinlan} and \cite{dSPlargerank}, two independent proofs are given of the fact that every linear subspace of $\Mat_n(\K)$ with a trivial spectrum has a dimension lesser than or equal to $\binom{n}{2}$. In \cite{dSPaffnonsingular}, spaces with a trivial spectrum whose dimension reaches this upper bound are classified up to similarity, extending Gerstenhaber's theorem (provided $\card \K>2$).

\paragraph{}
Here are our main results:

\begin{theo}\label{majotheo}
Let $V$ be a $1$-spec subspace of $\Mat_n(\K)$. \\
If $\car(\K)=2$ and $n=2$, then $\dim V \leq 3$. Otherwise $\dim V \leq 1+\binom{n}{2}$.
\end{theo}

\noindent The following corollary is trivial but must be stated:

\begin{cor}\label{majocor}
Let $V$ be a $\overline{1}$-spec subspace of $\Mat_n(\K)$. \\
If $\car(\K)=2$ and $n=2$, then $\dim V \leq 3$. Otherwise $\dim V \leq 1+\binom{n}{2}$.
\end{cor}

Note that if $\car(\K) \neq 2$ or $n \neq 2$,
the subspace $\K I_n+\NT_n(\K)$ of all upper triangular matrices with equal diagonal entries
is a $\overline{1}$-spec subspace of dimension $1+\binom{n}{2}$. \\
Moreover, setting
$$\frak{sl}_n(\K):=\bigl\{M \in \Mat_n(\K) : \; \tr(M)=0\bigr\},$$
we see that if $\car(\K)=2$, then a matrix of $\Mat_2(\K)$ has exactly one eigenvalue
in $\Kbar$ if and only if its trace is zero: it follows that
$\frak{sl}_2(\K)$ is a $3$-dimensional
$\overline{1}$-spec subspace of $\Mat_2(\K)$. The above upper bounds are therefore tight, and
we have the following full description of the $\overline{1}$-spec subspaces of $\Mat_2(\K)$
when $\car(\K)=2$:

\begin{prop}\label{egaltheo2}
If $\car(\K)=2$, then the $\overline{1}$-spec linear subspaces of $\Mat_2(\K)$ are the linear subspaces of $\frak{sl}_2(\K)$.
\end{prop}

We turn to the description of the $\overline{1}$-spec subspaces of $\Mat_n(\K)$ with maximal dimension.
Here is the most general situation:

\begin{theo}\label{egaltheo}
Let $V$ be a $\overline{1}$-spec subspace of $\Mat_n(\K)$
such that $\dim V=1+\binom{n}{2}$, with $n \geq 3$.
If $n \geq 5$ or $\car(\K) \wedge n=1$, then $V \sim \K I_n+\NT_n(\K)$.
\end{theo}

\begin{theo}\label{egaltheo4}
Assume $\car(\K)=2$. Then, up to similarity, there are exactly two
$\overline{1}$-spec subspaces of $\Mat_4(\K)$ with dimension $7$:
$\K I_4+\NT_4(\K)$ and the linear subspace
$$\calH:=\K I_4+\Biggl\{\begin{bmatrix}
0 & l_1 & l_2 & x \\
c_2 & 0 & y & c_1 \\
c_1 & x & 0 & c_2 \\
y & l_2 & l_1 & 0
\end{bmatrix} \mid (l_1,l_2,c_1,c_2,x,y)\in \K^6\Biggr\}.$$
\end{theo}

\begin{theo}
Assume $\K$ is an algebraically closed field of characteristic $3$.
Then, up to similarity, there are exactly three $\overline{1}$-spec subspaces of $\Mat_3(\K)$ with dimension $4$:
$\K I_3+\NT_3(\K)$ and the two subspaces
$$\biggl\{
\begin{bmatrix}
t & x & 0 \\
0 & t & y \\
z & 0 & t
\end{bmatrix} \mid (t,x,y,z)\in \K^4\biggr\} \quad \text{and} \quad
\biggl\{
\begin{bmatrix}
t & x & z \\
-z & t & y \\
x & 0 & t
\end{bmatrix} \mid (t,x,y,z)\in \K^4\biggr\}.$$
\end{theo}

For an arbitrary field of characteristic $3$, the case $n=3$ is far more complicated:
we wait until Section \ref{sectioncar3} to state the precise results.

We do not know yet how to classify the $1$-spec subspaces of $\Mat_n(\K)$ with maximal dimension,
although the following conjecture seems reasonable and would, if true, solve the question when
$\card \K>2$ and $n>4$ (using the results of \cite{dSPaffnonsingular}):

\begin{conj}
Let $V$ be a $1$-spec subspace of $\Mat_n(\K)$ such that $\dim V=1+\binom{n}{2}$.
Assume $n>4$. Then there exists a linear subspace $W$ of $\Mat_n(\K)$ with a trivial spectrum such that
$V=\K I_n+W$.
\end{conj}

\begin{proof}[Proof of Corollary \ref{majocor} and Theorem \ref{egaltheo} in the case $\car(\K) \wedge n=1$] ${}$ \\
Assume $\car(\K) \wedge n=1$. Then the results from Corollary \ref{majocor} and Theorem \ref{egaltheo}
are easy consequences of Gerstenhaber's theorem. Let indeed $V$ be a $\overline{1}$-spec subspace of $\Mat_n(\K)$.
Then $W:=\Ker(\tr_{|V})$ is a nilpotent subspace of $\Mat_n(\K)$ (since $\car(\K) \wedge n=1$) and $\codim_V W \leq 1$.
Gerstenhaber's theorem shows that $\dim W \leq \binom{n}{2}$, and hence
$$\dim V=\codim_V W+\dim W \leq 1+\binom{n}{2}.$$
Assume now that $\dim V=1+\binom{n}{2}$. Then $W$ is a hyperplane of $V$. Therefore $\dim W=\binom{n}{2}$ and
$W \sim \NT_n(\K)$ (see \cite{Serezhkin}).
Note also that the above inequality shows that $V$ is maximal among the $\overline{1}$-spec subspaces of $\Mat_n(\K)$. However $\K I_n+V$ is a $\overline{1}$-spec subspace
containing $V$, which shows that $I_n \in V$. Since $\car(\K) \wedge n=1$, one has
$I_n \not\in W$ hence $V=\K I_n+W$, and we deduce that $V \sim \K I_n+\NT_n(\K)$.
\end{proof}

\paragraph{Structure of the article:} ${}$ \\
The article has two main parts. In Section \ref{upperboundsection}, we prove
Theorem \ref{majotheo}, using arguments that are very similar to the ones used in
Theorem 9 of \cite{dSPlargerank}. \\
The remaining sections deal with the classification of $\overline{1}$-spec subspaces with maximal dimension.
Section \ref{maxdimensionsection} is devoted to a proof of Theorems \ref{egaltheo} and \ref{egaltheo4}:
adapting some ideas of \cite{dSPaffnonsingular}, we combine a key lemma (Proposition \ref{linelemma})
from Section \ref{upperboundsection} with Gerstenhaber's theorem in order to sort out the structure
of $\overline{1}$-spec subspaces with maximal dimension. We first work out the case $n \geq 5$,
and then the case $n=4$ and $\car(\K)=2$. In Section \ref{sectioncar3}, we use
the same line of reasoning to solve the case $n=3$ and $\car(\K)=3$: in that one,
many exceptional solutions arise, and the main difficulty lies in determining
necessary and sufficient conditions for two of them to be similar.
As we shall see, the classification depends on some arithmetic properties of the field $\K$.

\section{An upper bound for the dimension of a $1$-spec subspace}\label{upperboundsection}

\subsection{On the rank $1$ matrices in a $1$-spec subspace}

\begin{Not}
Let $V$ be a linear subspace of $\Mat_n(\K)$, and let $X \in \K^n \setminus \{0\}$. \\
We set $V_X:=\bigl\{M \in V : \; \im M \subset \K X\bigr\}$.
\end{Not}

Our proofs of Theorems \ref{majotheo} and \ref{egaltheo} are based on the following result:

\begin{prop}\label{linelemma}
Let $V$ be a $1$-spec subspace of $\Mat_n(\K)$, with $n \geq 2$.
Then:
\begin{enumerate}[(i)]
\item either $n=2$, $\car(\K)=2$ and $V=\frak{sl}_2(\K)$;
\item or there exists $X \in \K^n \setminus \{0\}$ such that $V_X=\{0\}$.
\end{enumerate}
\end{prop}

The proof involves the following result from \cite{dSPlargerank} (Proposition 10):

\begin{lemme}\label{lemmanontrivial}
Let $V$ be a linear subspace of $\Mat_n(\K)$ with a trivial spectrum.
Then there exists $i \in \lcro 1,n\rcro$ such that $V_{e_i}=\{0\}$, where $(e_1,\dots,e_n)$
denotes the canonical basis of $\K^n$.
\end{lemme}

\noindent Let us prove a corollary of it:

\begin{cor}\label{linecor}
Let $V$ be a linear subspace of $\Mat_n(\K)$ with a trivial spectrum.
Then there exists a basis $(f_1,\dots,f_n)$ of $\K^n$ such that
$V_{f_i}=\{0\}$ for every $i \in \lcro 1,n\rcro$.
\end{cor}

\begin{proof}
Denote by $F$ the linear subspace of $\K^n$ spanned by the non-zero vectors $X$ such that $V_X=\{0\}$. \\
Applying Lemma \ref{lemmanontrivial} to all the subspaces that are similar to $V$
shows that every basis of $\K^n$ contains a vector of $F$.
Classically, this shows that $F=\K^n$, which proves the claimed result: indeed, if the contrary holds,
then $F$ is included in a linear hyperplane $H$ of $E$; since $\GL_n(\K)$
acts transitively on the set of linear hyperplanes of $E$, we may find a basis $(f_1,\dots,f_n)$ of
$\K^n$ such that $H$ is defined by the equation $x_1+\cdots+x_n=0$ in this basis,
hence none of the vectors $f_1,\dots,f_n$ belongs to $F$, a contradiction.
\end{proof}

We may now prove Proposition \ref{linelemma}.

\begin{proof}[Proof of Proposition \ref{linelemma}]
We assume that $V_X \neq \{0\}$ for every $X \in \K^n \setminus \{0\}$, and prove that (i) holds. \\
Denote by $(e_1,\dots,e_{n-1})$ the canonical basis of $\K^{n-1}$. We naturally identify $\K^{n-1}$ with the subspace
$\K^{n-1} \times \{0\}$ of $\K^n$. \\
For $(i,j)\in \lcro 1,n\rcro^2$, denote by $E_{i,j}$ the elementary matrix of $\Mat_n(\K)$
with entry $1$ at the $(i,j)$-spot and zeroes elsewhere.
Denote by $W$ the linear subspace of $V$ consisting of its matrices with zero as the last row,
and note that $0$ is an eigenvalue of every matrix of $W$, hence $W$ has a trivial spectrum.
For $M \in W$, write
$$M=\begin{bmatrix}
K(M) & ? \\
0 & 0
\end{bmatrix} \quad \text{with $K(M) \in \Mat_{n-1}(\K)$.}$$
Then $K(W)$ is a linear subspace of $\Mat_{n-1}(\K)$ with a trivial spectrum.
Applying Corollary \ref{linecor} to $K(W)$, we find $P \in \GL_{n-1}(\K)$ such that, for
$F:=P K(W) P^{-1}$, one has $F_{e_i}=\{0\}$ for every $i \in \lcro 1,n-1\rcro$. \\
Set $V':=Q V Q^{-1}$, where $Q:=\begin{bmatrix}
P & 0 \\
0 & 1
\end{bmatrix}$. Then the assumptions show that $V'_X \neq \{0\}$ for every $X \in \K^n \setminus \{0\}$. \\
Let $i \in \lcro 1,n-1\rcro$. Then $V'_{e_i}$ contains a non-zero matrix $M$, which must belong to $QWQ^{-1}$:
however $F_{e_i}=\{0\}$ hence the first $n-1$ columns of $M$ are zero, which shows that $M$
is a scalar multiple of $E_{i,n}$. It follows that $E_{i,n} \in V'$ for every $i \in \lcro 1,n-1\rcro$,
hence $V'$ contains $\Vect(E_{1,n},\dots,E_{n-1,n})$.
Therefore $V$ contains $\Vect(E_{1,n},\dots,E_{n-1,n})$.
Conjugating $V$ with arbitrary matrices of $\GL_n(\K)$, this generalizes as follows:
$V$ contains every rank $1$ matrix with zero trace. \\
Since $V$ is a linear subspace, it must then contain the matrix $A:=E_{1,2}+E_{2,1}$,
which has $(x^2-1)x^{n-2}$ as characteristic polynomial.
This polynomial must have at most one root in $\K$, which shows that $n=2$ and $\car(\K)=2$.
Thus $V$ contains the basis $(E_{1,2},E_{2,1},J)$ of $\frak{sl}_2(\K)$,
where $J:=\begin{bmatrix}
1 & 1 \\
1 & 1
\end{bmatrix}$, hence $\frak{sl}_2(\K) \subset V \subset \Mat_2(\K)$, the last
inclusion being sharp since obviously $E_{1,1} \not\in V$.
Therefore $V=\frak{sl}_2(\K)$, which concludes our proof.
\end{proof}

\subsection{Proof of Theorem \ref{majotheo}}\label{proofoftheo1}

Let $V$ be a $1$-spec subspace of $\Mat_n(\K)$. \\
If $n=1$, we trivially have $\dim V \leq 1=1+\binom{1}{2}$. \\
If $n=2$ and $\car(\K) \neq 2$, then $V \subsetneq \Mat_2(\K)$ since $V$ does not contain $\Diag(1,0)$, and hence $\dim V \leq 3$. \\
Assume that $n \geq 3$ or that $n=2$ and $\car(\K)\neq 2$.
Then Proposition \ref{linelemma} shows that $V_X=\{0\}$ for some $X \in \K^n \setminus \{0\}$.
Conjugating $V$ with a well-chosen non-singular matrix, we lose no generality in assuming that $V_{e_n}=\{0\}$,
where $e_n=\begin{bmatrix}
0 & \cdots & 0 & 1
\end{bmatrix}^T$. \\
For $M \in V$, denote by $C(M)$ its last column, and set $U:=\{M \in V : \; C(M)=0\}$.
For $M \in U$, write
$$M=\begin{bmatrix}
J(M) & 0 \\
? & 0
\end{bmatrix} \quad \text{with $J(M) \in \Mat_{n-1}(\K)$.}$$
Since $V_{e_n}=\{0\}$, the linear map $M \mapsto J(M)$ is one-to-one on $U$, hence the rank theorem shows that
$$\dim V \leq n+\dim U=n+\dim J(U).$$
However $V$ is a $1$-spec subspace of $\Mat_n(\K)$, so, for every $M \in U$, the matrix $J(M)$ cannot have a non-zero eigenvalue in $\K$.
Theorem 9 of \cite{dSPlargerank} then shows that $\dim J(U) \leq \binom{n-1}{2}$.
Therefore $\dim V \leq n+\binom{n-1}{2}=1+\binom{n}{2}$, as claimed.

\begin{Rem}
Using a similar line of reasoning, it may be proven that $\frak{sl}_2(\K)$ is the sole
$3$-dimensional $1$-spec subspace of $\Mat_2(\K)$ if $\car(\K)=2$.
\end{Rem}

\section{On $\overline{1}$-spec subspaces with maximal dimension}\label{maxdimensionsection}

Let $n \geq 2$ and $V$ be a $\overline{1}$-spec subspace of $\Mat_n(\K)$ such that $\dim V=1+\binom{n}{2}$.
If $\car(\K)\wedge n=1$,
then we already know that Theorem \ref{egaltheo} holds for $V$ (see the end of Section \ref{pureintro}).
In the rest of the proof, we assume that $\car(\K)\mid n$.
The case $n=2$ has already been tackled (see Proposition \ref{egaltheo2}). \\
From now on, we assume that $n \geq 3$.

\paragraph{}
Notice that $\K I_n+V$ is a $\overline{1}$-spec subspace of $\Mat_n(\K)$
containing $V$, hence $V=\K I_n+V$ by Corollary \ref{majocor}, and in particular $I_n \in V$.
\paragraph{}
Finally, we will need the following notation and the subsequent remarks:

\begin{Not}
Given $M \in \Mat_n(\K)$, with characteristic polynomial $x^n+\underset{k=0}{\overset{n-1}{\sum}}a_k x^k \in \K[x]$,
we set $c_k(M):=(-1)^k a_{n-k}$ for every $k \in \lcro 1,n\rcro$.
\end{Not}

Classically, $c_2$ is a quadratic form on $\Mat_n(\K)$ and its polar form $b_2$, defined as
$b_2(A,B):=c_2(A+B)-c_2(A)-c_2(B)$, satisfies:
$$\forall (A,B)\in \Mat_n(\K)^2, \; b_2(A,B)=\tr(A)\tr(B)-\tr(AB).$$
Since $\car(\K) \mid n$, every matrix of $V$ has trace $0$, therefore
$$\forall (A,B)\in V^2, \; b_2(A,B)=-\tr(AB).$$
Notice that every singular matrix of $V$ is automatically nilpotent, which leads to the following result:

\begin{lemme}\label{singlemma}
Let $(A,B)\in V^2$. Assume that $A$, $B$ and $A+B$ are singular. \\
Then $\tr(AB)=0$.
\end{lemme}

\subsection{Setting things up}\label{setup}

We start with the same line of reasoning as in Section \ref{majotheo}. Since $n \geq 3$,
Proposition \ref{linelemma} shows that we may replace $V$ with a similar subspace so as to have, for
$e_n=\begin{bmatrix}
0 & \cdots & 0 & 1
\end{bmatrix}^T$,
\begin{equation}\label{cond1}
V_{e_n}=\{0\}.
\end{equation}
For $M \in V$, denote by $C(M)$ its last column, and set $Z:=\{M \in V : \; C(M)=0\}$.
For $M \in Z$, write
$$M=\begin{bmatrix}
J(M) & 0 \\
? & 0
\end{bmatrix} \quad \text{with $J(M) \in \Mat_{n-1}(\K)$.}$$
Note that $J(Z)$ is a nilpotent subspace of $\Mat_{n-1}(\K)$.
Since $V_{e_n}=\{0\}$, the map $J$ is one-to-one.
The rank theorem shows that $\dim V=\dim C(V)+\dim Z=\dim C(V)+\dim J(Z)$. However $\dim C(V) \leq n$ and $\dim J(Z) \leq \binom{n-1}{2}$
whilst $\dim V=n+\binom{n-1}{2}$. Therefore
$$\dim C(V)=n \quad \text{and} \quad \dim J(Z)=\binom{n-1}{2}.$$
Applying Gerstenhaber's theorem to $J(Z)$, we obtain $P \in \GL_{n-1}(\K)$ such that
$PJ(Z)P^{-1}=\NT_{n-1}(\K)$.
Set then $P_1:=\begin{bmatrix}
P & 0 \\
0 & 1
\end{bmatrix} \in \GL_n(\K)$ and replace $V$ with $P_1 VP_1^{-1}$.
Notice that $P_1^{-1}e_n=e_n$, therefore condition \eqref{cond1} is still satisfied in this new setting,
but we now have
\begin{equation}\label{cond2}
J(Z)=\NT_{n-1}(\K).
\end{equation}

\paragraph{}
Denote now by $W$ the set of all matrices of $V$ with zero as first row.
For $M \in W$, write
$$M=\begin{bmatrix}
0 & 0 \\
? & R(M)
\end{bmatrix} \quad \text{with $R(M) \in \Mat_{n-1}(\K)$.}$$
Note the following obvious properties of $R(W)$:
\begin{enumerate}[(i)]
\item $R(W)$ is a nilpotent subspace of $\Mat_{n-1}(\K)$;
\item The shape of $J(Z)$ shows that, for every $N \in \NT_{n-2}(\K)$, the subspace
$R(W)$ possesses a matrix of the form
$\begin{bmatrix}
N & 0 \\
? & 0
\end{bmatrix}$.
\end{enumerate}

Let $C_1 \in \Mat_{n-2,1}(\K)$. Since $\dim C(V)=n$, we know that $V$ contains a matrix of the form
$\begin{bmatrix}
? & ? & 0 \\
? & ? & C_1 \\
? & ? & 0
\end{bmatrix}$.
By adding a well-chosen matrix of $Z$ and a scalar multiple of $I_n$, we deduce that $W$ contains a matrix of the form
$\begin{bmatrix}
0 & 0 & 0 \\
? & ? & C_1 \\
? & ? & ?
\end{bmatrix}$, i.e., $R(W)$ contains a matrix of the form
$\begin{bmatrix}
? & C_1 \\
? & ?
\end{bmatrix}$. Since this holds for every $C_1 \in \Mat_{n-2,1}(\K)$, combining this with
point (ii) above yields that $\dim R(W) \geq (n-2)+\binom{n-2}{2}=\binom{n-1}{2}$,
and hence $\dim R(W)=\binom{n-1}{2}$ by Gerstenhaber's theorem.
Since $R(W)$ is nilpotent, no matrix of it has $\begin{bmatrix}
0 & \cdots & 0 & 1
\end{bmatrix}^T$ as the last column. By the factorization lemma for linear maps, we deduce that there
exists a (unique) $L_1 \in \Mat_{1,n-2}(\K)$ such that every matrix $M$ of $R(W)$ has $\begin{bmatrix}
C_M \\
L_1C_M
\end{bmatrix}$ as the last column (for some $C_M \in \Mat_{n-2,1}(\K)$).

\begin{claim}
For $Q:=\begin{bmatrix}
I_{n-2} & 0 \\
-L_1 & 1
\end{bmatrix}$, one has $QR(W)Q^{-1}=\NT_{n-1}(\K)$.
\end{claim}

\begin{proof}
Set $\calT:=QR(W)Q^{-1}$. With the above results on $R(W)$, we find:
\begin{itemize}
\item[(a)] Every matrix of $\calT$ has entry $0$ at the $(n-1,n-1)$-spot and, for every
$C_1 \in \Mat_{n-2,1}(\K)$, the subspace $\calT$ contains a matrix with
 $\begin{bmatrix}
C_1 \\
0
\end{bmatrix}$ as the last column;
\item[(b)] For every $N \in \NT_{n-2}(\K)$, the subspace $\calT$ contains a matrix of the form
$\begin{bmatrix}
N & 0 \\
? & 0
\end{bmatrix}$.
\end{itemize}
Moreover, $\calT$ is a $\binom{n-1}{2}$-dimensional nilpotent subspace of $\Mat_{n-1}(\K)$.
By Gerstenhaber's theorem, there is an increasing sequence $\{0\}=E_0 \subsetneq E_1 \subsetneq \cdots \subsetneq E_{n-1}=\K^{n-1}$
of linear subspaces such that
$\forall k \in \lcro 1,n-1\rcro, \; \forall M \in \calT, \; M(E_k) \subset E_{k-1}$. \\
Point (a) then yields $E_{n-2}=\K^{n-2}\times \{0\}$. It follows that every matrix of $\calT$ has
$0$ as the last row, hence point (b) may be refined as follows:
\begin{itemize}
\item[(b')] For every $N \in \NT_{n-2}(\K)$, the subspace $\calT$ contains
$\widetilde{N}:=\begin{bmatrix}
N & 0 \\
0 & 0
\end{bmatrix}\in \Mat_{n-1}(\K)$.
\end{itemize}
However $\K^{n-3} \times \{0\} \subset \underset{N \in \NT_{n-2}(\K)}{\sum} \widetilde{N}(E_{n-2}) \subset E_{n-3}$, and $\dim(\K^{n-3} \times \{0\})=\dim E_{n-3}$,
therefore $E_{n-3}=\K^{n-3} \times \{0\}$.
Continuing by downward induction, we find that $E_k=\K^k \times \{0\}$ for every $k \in \lcro 0,n-1\rcro$,
which shows that $\calT \subset \NT_{n-1}(\K)$, and hence $\calT=\NT_{n-1}(\K)$
since $\dim \calT=\binom{n-1}{2}=\dim \NT_{n-1}(\K)$.
\end{proof}

We now set $Q_1:=\begin{bmatrix}
1 & 0 \\
0 & Q
\end{bmatrix} \in \GL_n(\K)$ and replace $V$ with $Q_1 V Q_1^{-1}$:
again, condition \eqref{cond1} still holds in this new setting, and $J(Z)$ has not been modified, therefore
condition \eqref{cond2} still holds. Since $Q_1$ stabilizes $\{0\} \times \K^{n-1}$, the subspaces
$W$ and $R(W)$ have been replaced respectively with
$Q_1 W Q_1^{-1}$ and $Q R(W) Q^{-1}$, therefore we now have:
\begin{equation}\label{cond3}
R(W)=\NT_{n-1}(\K).
\end{equation}

\subsection{Special matrices in $V$}\label{specialmat}

Using conditions \eqref{cond1} and \eqref{cond2}, we find that, for every $N \in \NT_{n-1}(\K)$,
the subspace $V$ contains a unique matrix of the form
$\begin{bmatrix}
N & 0 \\
? & 0
\end{bmatrix}$.
We deduce that:
\begin{itemize}
\item There are two linear maps $\varphi : \Mat_{1,n-2}(\K) \rightarrow \Mat_{1,n-2}(\K)$
and $f : \Mat_{1,n-2}(\K) \rightarrow \K$ such that, for every $L \in \Mat_{1,n-2}(\K)$,
the subspace $V$ contains the matrix
$$A_L:=\begin{bmatrix}
0 & L & 0 \\
0 & 0_{n-2} & 0 \\
f(L) & \varphi(L) & 0
\end{bmatrix};$$
\item There is a linear form $h : \NT_{n-2}(\K) \rightarrow \K$ such that, for every $U \in \NT_{n-2}(\K)$,
the subspace $V$ contains the matrix
$$E_U:=\begin{bmatrix}
0 & 0 & 0 \\
0 & U & 0 \\
h(U) & 0 & 0
\end{bmatrix}$$
(there, we have also used condition \eqref{cond3}).
\end{itemize}
\vskip 3mm
Since $\dim C(V)=n$, we know that some matrix of $V$ has entry $1$ at the $(1,n)$-spot.
By linearly combining such a matrix with $I_n$ and a well-chosen $A_L$, we deduce that
every row matrix $L' \in \Mat_{1,n}(\K)$ is the first row of some matrix of $V$.
Denote by $G$ the linear subspace of $V$ consisting of the matrices $M \in V$ with $\tr M=0$ and all columns zero starting from the second one.
Applying the rank theorem then shows that
$\dim V=n+\dim R(W)+\dim G$.
However $n+\dim R(W)=n+\binom{n-1}{2}=\dim V$, therefore $G=\{0\}$. \\
Condition \eqref{cond3} then yields
that, for every $N \in \NT_{n-1}(\K)$, the linear subspace $V$ contains a unique matrix of the form
$\begin{bmatrix}
0 & 0 \\
? & N
\end{bmatrix}$.
We deduce a new family of matrices in $V$:
\begin{itemize}
\item There are two linear maps $\psi : \Mat_{n-2,1}(\K) \rightarrow \Mat_{n-2,1}(\K)$
and $g : \Mat_{n-2,1}(\K) \rightarrow \K$ such that, for every $C \in \Mat_{n-2,1}(\K)$,
the subspace $V$ contains the matrix
$$B_C:=\begin{bmatrix}
0 & 0 & 0 \\
\psi(C) & 0_{n-2} & C \\
g(C) & 0 & 0
\end{bmatrix}.$$
\end{itemize}
Finally, we have seen that some matrix of $V$ has $\begin{bmatrix}
0 & \cdots & 0 & 1
\end{bmatrix}$ as first row. By adding to it a well-chosen matrix of the form $B_C$, we find a matrix
$$J=\begin{bmatrix}
0 & 0_{1,n-2} & 1 \\
? & ? & 0_{n-2,1} \\
? & ? & ?
\end{bmatrix} \in V.$$

\vskip 3mm
Obviously the subspace $\{A_L \mid L \in \Mat_{1,n-2}(\K)\}+\{B_C \mid C \in \Mat_{n-2,1}(\K)\}
+\{E_U \mid U \in \NT_{n-2}(\K)\}+\Vect(J,I_n)$ of $V$ has dimension $2(n-2)+\binom{n-2}{2}+2=\dim V$,
hence $V$ is spanned by the above matrices.
At this point, we need to examine three cases separately: $n\geq 5$, $n=4$ and $n=3$
(the last one is dealt with in Section \ref{sectioncar3}).

\subsection{The case $n \geq 5$}

\subsubsection{Analyzing $\varphi$ and $\psi$}\label{phietpsi}

\begin{claim}
There exists $\lambda \in \K$ such that $\varphi=\lambda\,\id$ and $\psi=-\lambda\,\id$.
\end{claim}

Notice first that given $(L,C)\in \Mat_{1,n-2}(\K) \times \Mat_{n-2,1}(\K)$, one has $\rk(A_L) \leq 2$
and $\rk(B_C) \leq 2$, and hence $\rk(A_L+B_C) \leq 4$: the matrix $A_L+B_C$, which belongs to $V$,
is singular and therefore nilpotent.

\begin{proof}
Let $(L,C)\in \Mat_{1,n-2}(\K) \times \Mat_{n-2,1}(\K)$. \\
Assume that $LC=0$. Then
$$(A_L+B_C)^2=\begin{bmatrix}
L\psi(C) & 0 & 0 \\
? & \psi(C)L+C\varphi(L) & 0 \\
? & ? & \varphi(L)C
\end{bmatrix}.$$
Since $A_L+B_C$ is nilpotent, we deduce that $L\psi(C)=\varphi(L)C=0$. \\
We equip $\Mat_{1,n-2}(\K)$ with the non-degenerate symmetric bilinear form $(L_1,L_2) \mapsto L_1L_2^T$. \\
Let $L \in \Mat_{1,n-2}(\K)$. Then the above result yields $\varphi(L) \in (\{L\}^\bot)^\bot=\Vect(L)$. \\
It follows that every non-zero vector of $\Mat_{1,n-2}(\K)$ is an eigenvector of $\varphi$,
which classically yields that $\varphi=\lambda\, \id$ for some $\lambda \in \K$. \\
With the same line of reasoning, we find that $\psi=\mu\, \id$ for some $\mu \in \K$. \\
Choose $(L,C)\in \Mat_{1,n-2}(\K) \times \Mat_{n-2,1}(\K)$ such that $LC \neq 0$. Then $A_L$, $B_C$ and $A_L+B_C$ are all singular, therefore $\tr(A_LB_C)=0$ by Lemma \ref{singlemma}, i.e., $\psi(L)C+L \varphi(C)=0$.
We deduce that $(\lambda+\mu)\,LC=0$, and hence $\mu=-\lambda$.
\end{proof}

\subsubsection{One last conjugation}\label{lastconjugate}

Set $P':=\begin{bmatrix}
1 & 0 & 0 \\
0 & I_{n-2} & 0 \\
\lambda & 0 & 1
\end{bmatrix} \in \GL_n(\K)$, and note that, for every $(L,C,U) \in \Mat_{1,n-2}(\K) \times \Mat_{n-2,1}(\K) \times \NT_{n-2}(\K)$:
$$(P')^{-1}A_L P'=\begin{bmatrix}
0 & L & 0 \\
0 & 0_{n-2} & 0 \\
f(L) & 0 & 0
\end{bmatrix} \quad ; \quad (P')^{-1}B_C P'=\begin{bmatrix}
0 & 0 & 0 \\
0 & 0_{n-2} & C \\
g(C) & 0 & 0
\end{bmatrix}$$
and
$$(P')^{-1}E_U P'=\begin{bmatrix}
0 & 0 & 0 \\
0 & U & 0 \\
h(U) & 0 & 0
\end{bmatrix}.$$

Notice that $P'e_n=e_n$ hence $(P')^{-1}V P'$ still satisfies condition \eqref{cond1}.
The above matrices show that condition \eqref{cond2} and \eqref{cond3} are obviously satisfied.
Replacing $V$ with $(P')^{-1}V P'$, we thus preserve all the previous conditions but we now have
the additional one:
$$\lambda=0.$$
At this point, our aim is to prove that $V \subset \K I_n+\NT_n(\K)$,
which will suffice since $V$ and $\K I_n+\NT_n(\K)$ have the same dimension.
In order to do so, we prove that every matrix of the type $A_L$, $B_C$, $E_U$ or $J$ is strictly upper-triangular:
this suffices to prove our claim since these matrices, together with $I_n$, span $V$.

\subsubsection{Analyzing $f$, $g$ and $h$}

\begin{claim}\label{f=g=0}
One has $f=0$ and $g=0$.
\end{claim}

\begin{proof}
Let $(L,C)\in \Mat_{1,n-2}(\K) \times \Mat_{n-2,1}(\K)$. Set $\alpha:=f(L)+g(C)$.
Then a straightforward computation shows that
$$(A_L+B_C)^3=\begin{bmatrix}
\alpha\,LC & 0 & 0 \\
0 & \alpha\,CL & 0 \\
0 & 0 & \alpha\,LC
\end{bmatrix}.$$
However $A_L+B_C$ is nilpotent (see Paragraph \ref{phietpsi}), therefore $\alpha\,LC=0$.
It follows that $LC \neq 0 \Rightarrow f(L)+g(C)=0$. \\
Let $C_1 \in \Mat_{1,n-2}(\K)$. Since $n-2 \geq 2$, we may choose
$L \in \Mat_{1,n-2}(\K) \setminus \{0\}$ such that $LC_1=0$.
Then the linear map $g$ is constant on the affine hyperplane $\bigl\{C \in \Mat_{1,n-2}(\K): \; LC=1\bigr\}$,
therefore $g$ vanishes on its translation vector space $\bigl\{C \in \Mat_{1,n-2}(\K): \; LC=0\bigr\}$.
In particular $g(C_1)=0$. We deduce that $g=0$. The same line of reasoning yields $f=0$.
\end{proof}

\begin{claim}
One has $h=0$.
\end{claim}

\begin{proof}
Let $U \in \NT_{n-2}(\K)$ such that $\rk U=1$. Set $\beta:=h(U)$. Note that $U^2=0$ since $U$ is a rank $1$ nilpotent matrix.
Let $(L,C) \in \Mat_{1,n-2}(\K) \times \Mat_{n-2,1}(\K)$.
Set
$$M:=A_L+B_C+E_U=\begin{bmatrix}
0 & L & 0 \\
0 & U & C \\
\beta & 0 & 0
\end{bmatrix}.$$
A straightforward computation (using $U^2=0$) yields
$$M^3=\begin{bmatrix}
\beta\, LC & 0 & LUC \\
\beta\, UC & \beta\, CL & 0 \\
0 & \beta\,LU & \beta\, LC
\end{bmatrix}.$$
Note that $\rk A_L \leq 1$, $\rk B_C \leq 1$ and $\rk E_U \leq 2$, therefore
$\rk M \leq 4$. Since $n \geq 5$, we deduce that $M$ is nilpotent.
It follows that $LUC=0 \Rightarrow \beta\,LC=0$. \\
Choosing $L:=\begin{bmatrix}
0 & \cdots & 0 & 1
\end{bmatrix}$ and $C:=L^T$, we then have $LUC=0$ (since $U$ is strictly upper-triangular)
whilst $LC=1$, and we deduce that $\beta=0$. \\
We have just established that the linear form $h$ vanishes on every rank $1$ matrix of $\NT_{n-2}(\K)$, which proves
our claim since $\NT_{n-2}(\K)$ is obviously spanned by its rank $1$ matrices.
\end{proof}

We thus have, for every $(L,C,U) \in \Mat_{1,n-2}(\K) \times \Mat_{n-2,1}(\K) \times \NT_{n-2}(\K)$,
$$A_L=\begin{bmatrix}
0 & L & 0 \\
0 & 0_{n-2} & 0 \\
0 & 0 & 0
\end{bmatrix} \; ; \; B_C=\begin{bmatrix}
0 & 0 & 0 \\
0 & 0_{n-2} & C \\
0 & 0 & 0
\end{bmatrix} \; \text{and} \; E_U=\begin{bmatrix}
0 & 0 & 0 \\
0 & U & 0 \\
0 & 0 & 0
\end{bmatrix}.$$
It now suffices to show that $J$ is strictly upper-triangular.

\subsubsection{Dissecting $J$}

Here, we aim at proving that the matrix $J$ of Section \ref{specialmat} may be chosen to have zero entries everywhere except at the $(1,n)$-spot.

\vskip 2mm
Adding to $J$ a well-chosen matrix of type $E_U$, we may assume
that
$$J=\begin{bmatrix}
0 & 0 & 1 \\
C_1 & T & 0 \\
b & L_1 & a
\end{bmatrix} \quad \text{where $L_1 \in \Mat_{1,n-2}(\K)$, $C_1 \in \Mat_{n-2,1}(\K)$, $(a,b)\in \K^2$}$$
and $T=(t_{i,j})$ is a \emph{lower-triangular} matrix of $\Mat_{n-2}(\K)$.

\begin{claim}
One has $C_1=0$ and $L_1=0$, whilst $T=\alpha\,I_{n-2}$ for some $\alpha \in \K$.
\end{claim}

\begin{proof}
Denote by $l$ the last entry of $L_1$.
Setting $L_2:=\begin{bmatrix}
0 & \cdots & 0 & 1
\end{bmatrix}$, we remark that
$V$ contains the matrix $M:=J+A_{L_2}+(l-a).I_n+(t_{n-2,n-2}+l-a)B_{L_2^T}$, which has
identical $(n-1)$-th and $n$-th columns (and the same non-diagonal entries as $J$ on the first $n-2$ rows)
and is therefore singular.
Choose $U \in \NT_{n-2}(\K)$ which has a last column equal to zero. Again, $M+E_U$ is singular (it has the same last two columns as $M$).
So are $E_U$ and $M$, therefore Lemma \ref{singlemma} yields $\tr(ME_U)=0$. \\
By varying $U$, we deduce that $t_{i,j}=0$ for every $(i,j)\in \lcro 1,n-2\rcro^2$ such that $j<i \leq n-3$. \\
With the same line of reasoning, we find that $\tr(MA_L)=0$ for every $L \in \Mat_{1,n-2}(\K)$ with $0$ as the last entry,
which shows that $C_1$ has zero entries from the first to the $(n-3)$-th one. \\
With the same line of reasoning, but replacing the last two columns with the first two rows,
we find that $t_{i,j}=0$ for every $(i,j)\in \lcro 1,n-2\rcro^2$ such that $2 \leq j<i$,
and $L_1$ has zero entries starting from the second one.

Denote by $(e_1,\dots,e_n)$ the canonical basis of $\K^n$.
With the above results, we find that $e_3,\dots,e_{n-1}$ are eigenvectors of $J$ with respective eigenvalues
$t_{2,2},\dots,t_{n-2,n-2}$, therefore $t_{2,2}=\cdots=t_{n-2,n-2}$. \\
Moreover $e_2,\dots,e_{n-2}$ are eigenvectors of $J^T$ with respective eigenvalues
$t_{1,1},\dots,t_{n-3,n-3}$, therefore $t_{1,1}=\cdots=t_{n-3,n-3}$. Since $n-2 \geq 3$, it follows that
the $t_{i,i}$'s are all equal to some $\alpha \in \K$.

Let us now consider $J':=J-\alpha.I_n \in V$. Note that $\rk(J') \leq 3$ since all the columns of $J'$ from the third one to the $(n-1)$-th one are zero.
Let $L \in \Mat_{1,n-2}(\K)$. Then $\rk(J'+A_L) \leq 4<n$. The matrices $J'$, $A_L$ and $J'+A_L$
are all singular, therefore Lemma \ref{singlemma} shows that $\tr(J'A_L)=0$
i.e., $LC_1=0$. Since this holds for every $L \in \Mat_{1,n-2}(\K)$, we find that $C_1=0$. \\
With the same line of reasoning (but using the $B_C$'s instead), we find $L_1=0$. \\
Finally, with the same line of reasoning with the elementary matrix $U=E_{1,n-2} \in \NT_{n-2}(\K)$
with zero entries everywhere except at the $(1,n-2)$-spot where the entry is $1$, we find
$\tr(J'U)=0$ (note that $\rk(U)=1$), i.e., $t_{n-2,1}=0$. With the above results, this finally shows that
$T=\alpha\,I_{n-2}$.
\end{proof}

\begin{Rem}
If $\car(\K) \neq 2$ (still assuming that $\car(\K)\mid n$), then the above proof may be greatly simplified. Indeed,
if a matrix of $\Mat_n(\K)$ has $\lambda \in \Kbar$ as sole eigenvalue, then its characteristic polynomial is
$\bigl((x-\lambda)^{n/p}\bigr)^p$, where $p:=\car(\K)$,
hence $c_2(M)=0$. It follows that $c_2$ vanishes everywhere on $V$, and therefore $\tr(MN)=0$ for every $(M,N)\in V^2$.
Applying this to $M=J$ and $N$ being of any one of the types $A_L$, $B_C$ and $E_U$, we find that
$C_1=0$, $L_1=0$ and $T$ is diagonal. We let the reader finish the proof in that case.
\end{Rem}

\begin{claim}
One has $a=b=\alpha=0$.
\end{claim}

\begin{proof}
Set $J':=J-\alpha.I_n$.
Then
$$J'=\begin{bmatrix}
-\alpha & 0 & 1 \\
0 & 0_{n-2} & 0 \\
b & 0 & a-\alpha
\end{bmatrix}$$
Then $J'$ is singular, and therefore nilpotent since $J' \in V$. Then $\tr(J')=0$ and $c_2(J')=0$, which shows that
$a-\alpha=\alpha$ and $b=-\alpha^2$.
Therefore
$$J'=\begin{bmatrix}
-\alpha & 0 & 1 \\
0 & 0_{n-2} & 0 \\
-\alpha^2 & 0 & \alpha
\end{bmatrix}.$$
Choose $(L,C)\in \Mat_{1,n-2}(\K) \times \Mat_{n-2,1}(\K)$ such that $LC=1$. \\
Set $M:=J'+A_L+B_C=\begin{bmatrix}
-\alpha & L & 1 \\
0 & 0_{n-2} & C \\
-\alpha^2 & 0 & \alpha
\end{bmatrix}$. Note that $\rk(M) \leq 3$. Hence $M$ is singular and therefore nilpotent.
A straightforward computation shows that the first column of $M^3$ is $\begin{bmatrix}
-\alpha^2 & 0 & \cdots & 0
\end{bmatrix}^T$, which yields $\alpha=0$. The conclusion easily follows since $a=2\alpha$ and $b=-\alpha^2$.
\end{proof}

Our proof is now complete: we know that $J$ is strictly upper-triangular,
and so is any matrix of type $A_L$, $B_C$ or $E_U$.
It follows that $V \subset \K I_n+\NT_n(\K)$, and the equality of spaces follows from the equality of their dimensions.
This finishes our proof of Theorem \ref{egaltheo}.

\subsection{The case $n=4$ and $\car(\K)=2$}

Recall from Theorem \ref{egaltheo4} the definition of
$$\calH=\K I_4+\Biggl\{\begin{bmatrix}
0 & l_1 & l_2 & x \\
c_2 & 0 & y & c_1 \\
c_1 & x & 0 & c_2 \\
y & l_2 & l_1 & 0
\end{bmatrix} \mid (l_1,l_2,c_1,c_2,x,y)\in \K^6\Biggr\},$$
which is obviously a $7$-dimensional linear subspace of $\Mat_4(\K)$.

\begin{claim}
The set $\calH$ is a $\overline{1}$-spec subspace of $\Mat_4(\K)$
and it is not similar to $\K I_4+\NT_4(\K)$.
\end{claim}

\begin{proof}
A straightforward computation shows that, for every $(l_1,l_2,c_1,c_2,x,y) \in \K^6$,
the matrix $\begin{bmatrix}
0 & l_1 & l_2 & x \\
c_2 & 0 & y & c_1 \\
c_1 & x & 0 & c_2 \\
y & l_2 & l_1 & 0
\end{bmatrix}$ has characteristic polynomial
$t^4+a$ (in $\K[t]$), with $a=\bigl((l_1+l_2)(c_1+c_2)+xy\bigr)^2+l_1c_2(x+y)^2$,
and has therefore a unique eigenvalue in $\Kbar$ (since $\car(\K)=2$).
It follows that $\calH$ is a $\overline{1}$-spec subspace. \\
If $\calH$ were similar to $\K I_4+\NT_4(\K)$, the set of singular matrices of $\calH$ would be a linear subspace of $\calH$.
However, $\begin{bmatrix}
0 & 0 & 0 & 1 \\
0 & 0 & 0 & 0 \\
0 & 1 & 0 & 0 \\
0 & 0 & 0 & 0
\end{bmatrix}$ and $\begin{bmatrix}
0 & 0 & 0 & 0 \\
0 & 0 & 1 & 0 \\
0 & 0 & 0 & 0 \\
1 & 0 & 0 & 0
\end{bmatrix}$ are singular, whereas their sum is not. Therefore $\calH$ is not similar to $\K I_4+\NT_4(\K)$.
\end{proof}

\paragraph{}
We now come right back to the end of Section \ref{specialmat} and try to prove that $V$
is similar to $\K I_4+\NT_4(\K)$ or to $\calH$. Notice first that, for every $M \in V$,
its characteristic polynomial has the form $(t+\lambda)^4=t^4+\lambda^4$ for some $\lambda \in \Kbar$,
therefore $c_2$ and $c_3$ vanish everywhere on $V$.

\begin{claim}
There is a (unique) matrix $A \in \Mat_2(\K)$ such that
$\forall L \in \Mat_{1,2}(\K), \; \varphi(L)=LA$ and $\forall C \in \Mat_{2,1}(\K), \; \psi(C)=AC$.
\end{claim}

\begin{proof}
Indeed, we know that there are two matrices $A$ and $B$ in $\Mat_2(\K)$ such that
$\forall L \in \Mat_{1,2}(\K), \; \varphi(L)=LA$ and $\forall C \in \Mat_{2,1}(\K), \; \psi(C)=BC$.
Since $c_2$ vanishes everywhere on $V$, we have $\tr(A_LB_C)=0$ for every $(L,C) \in \Mat_{1,2}(\K) \times \Mat_{2,1}(\K)$,
which yields $LAC+LBC=0$, i.e., $L(A+B)C=0$. It follows that $A+B=0$, hence $B=A$ since $\car(\K)=2$.
\end{proof}

\begin{claim}
One has $f=0$ and $g=0$.
\end{claim}

\begin{proof}
Let $(L,C)\in \Mat_{1,2}(\K) \times \Mat_{2,1}(\K)$. Set $\alpha:=f(L)+g(C)$ and
$$M:=A_L+B_C=\begin{bmatrix}
0 & L & 0 \\
AC & 0 & C \\
\alpha & LA & 0
\end{bmatrix}.$$
Notice that $c_3(M)=\tr(\com(M))$, where $\com(M)$ denotes the comatrix of $M$. \\
Obviously the first and last diagonal entries of $\com(M)$ are zero
hence
$$c_3(M)=\begin{vmatrix}
0 & l_2 & 0 \\
? & 0 & c_2 \\
\alpha & ? & 0
\end{vmatrix}+\begin{vmatrix}
0 & l_1 & 0 \\
? & 0 & c_1 \\
\alpha & ? & 0
\end{vmatrix}=\alpha\,LC,$$
where $L=\begin{bmatrix}
l_1 & l_2
\end{bmatrix}$ and $C=\begin{bmatrix}
c_1  \\
c_2
\end{bmatrix}$. Therefore $LC \neq 0 \Rightarrow f(L)+g(C)=0$.
As in the proof of Claim \ref{f=g=0}, we deduce that $f=0$ and $g=0$.
\end{proof}

Remark that $\NT_2(\K)$ is spanned by $K:=\begin{bmatrix}
0 & 1 \\
0 & 0
\end{bmatrix}$. Setting now $\alpha:=h(K)$ and $E:=\begin{bmatrix}
0 & 0 & 0 & 0 \\
0 & 0 & 1 & 0 \\
0 & 0 & 0 & 0 \\
\alpha & 0 & 0 & 0
\end{bmatrix}$, we deduce that
$V=\bigl\{A_L \mid L \in \Mat_{1,2}(\K)\}+\bigl\{B_C \mid C \in \Mat_{2,1}(\K)\}+\Vect(I_4,E,J)$.

\begin{claim}
There exists $\lambda \in \K$ such that
$A=\begin{bmatrix}
\lambda & \alpha \\
\alpha & \lambda
\end{bmatrix}$.
\end{claim}

\begin{proof}
Let $(L,C)\in  \Mat_{1,2}(\K) \times \Mat_{2,1}(\K)$ and set $M:=A_L+B_C+E$.
Then a straightforward computation shows that
$$c_3(M)=\alpha LC+l_1\,(AC)_2+c_2\,(LA)_1,$$
where $(AC)_2$ denotes the second entry of $AC$ and $(LA)_1$ the first one of $LA$. \\
Write $A=\begin{bmatrix}
a_{1,1} & a_{1,2} \\
a_{2,1} & a_{2,2}
\end{bmatrix}$. Notice then that $(L,C) \mapsto \alpha LC+l_1\,(AC)_2+c_2\,(LA)_1$ is a bilinear form whose matrix in the
respective canonical bases of $\Mat_{1,2}(\K)$ and $\Mat_{2,1}(\K)$ is
$$\alpha\,I_2+\begin{bmatrix}
a_{2,1} & a_{2,2} \\
0 & 0
\end{bmatrix}+\begin{bmatrix}
0 & a_{1,1} \\
0 & a_{1,2}
\end{bmatrix}=\begin{bmatrix}
\alpha+a_{2,1} & a_{1,1}+a_{2,2} \\
0 & \alpha+a_{1,2}
\end{bmatrix}.$$
Since $c_3$ vanishes everywhere on $V$, it follows that this matrix is zero, which yields $a_{2,1}=a_{1,2}=\alpha$ and $a_{2,2}=a_{1,1}$,
as claimed.
\end{proof}

\paragraph{}
As in Paragraph \ref{lastconjugate}, we now replace $V$ with $(P')^{-1} V P'$
where
$P':=\begin{bmatrix}
1 & 0 & 0 \\
0 & I_2 & 0 \\
\lambda & 0 & 1
\end{bmatrix} \in \GL_4(\K)$. Again, all the former conditions still hold in that case,
and we now have $A=\begin{bmatrix}
0 & \alpha \\
\alpha & 0
\end{bmatrix}$.

\paragraph{}
We finish by analyzing $J$.
By summing it with a well-chosen scalar multiple of $E$, we lose no generality in assuming that
$$J=\begin{bmatrix}
0 & 0 & 1 \\
C_1 & T & 0 \\
\beta & L_1 & a
\end{bmatrix}$$
where $(L_1,C_1)\in  \Mat_{1,2}(\K) \times \Mat_{2,1}(\K)$, $(a,\beta)\in \K^2$ and $T$ is a lower-triangular matrix of $\Mat_2(\K)$. \\
Since $c_2$ vanishes everywhere on $V$, we find that $\tr(JA_L)=\tr(JB_C)=\tr(JE)=0$ for every
$(L,C) \in \Mat_{1,2}(\K) \times \Mat_{2,1}(\K)$, which yields $C_1=0$, $L_1=0$ and $t_{2,1}=\alpha$.

\begin{claim}
One has $a=0$, and there exists $b \in \K$ such that $\beta=b^2$ and
$T=\begin{bmatrix}
b & 0 \\
\alpha & b
\end{bmatrix}$.
\end{claim}

\begin{proof}
Write $T=\begin{bmatrix}
b & 0 \\
\alpha & c
\end{bmatrix}$.
Denote by $(e_1,e_2,e_3,e_4)$ the canonical basis of $\K^4$. Note that the endomorphism $X \mapsto JX$ of $\K^4$ stabilizes both of the subspaces
$\Vect(e_2,e_3)$ and $\Vect(e_1,e_4)$ and the matrices of the induced endomorphisms in the respective bases $(e_2,e_3)$ and
$(e_1,e_4)$ are $T$ and $\begin{bmatrix}
0 & 1 \\
\beta & a
\end{bmatrix}$. Therefore those matrices have the same unique eigenvalue in $\overline{\K}$, which must be $b$.
This yields $b=c$, $a=0$ and $\beta=b^2$.
\end{proof}

\begin{claim}
One has $b=0$.
\end{claim}

\begin{proof}
Choose $(L,C)\in \Mat_{1,2}(\K) \times \Mat_{2,1}(\K)$ such that $LC \neq 0$.
A straightforward computation shows that
$$c_3(A_L+B_C+J)=b^2\,LC,$$
and hence $b=0$.
\end{proof}

We now have
$$V=\Biggl\{\begin{bmatrix}
\lambda & l_1 & l_2 & x \\
\alpha c_2 & \lambda & y & c_1 \\
\alpha c_1 & \alpha x & \lambda & c_2 \\
\alpha y & \alpha l_2 & \alpha l_1 & \lambda
\end{bmatrix} \mid (l_1,l_2,c_1,c_2,x,y,\lambda) \in \K^7\Biggr\}.$$
If $\alpha=0$, then we readily have $V=\K I_4+\NT_4(\K)$. \\
Assume finally that $\alpha \neq 0$. Then, for $D:=\Diag(1,1,1,\alpha)$,
a straightforward computation yields $D^{-1}\,V\,D=\calH$.
This completes the proof of Theorem \ref{egaltheo4}.

\section{On $\overline{1}$-spec subspaces of $\Mat_3(\K)$ when $\car(\K)=3$}\label{sectioncar3}

In this section, we assume that $\car(\K)=3$.

\subsection{Opening remarks}

Here, we will use considerations from Witt's theory of quadratic forms (see \cite[Chapters VII, VIII and IX]{invitquad}).

Since $\car(\K)=3$, a matrix $M \in \Mat_3(\K)$ has only one eigenvalue in $\Kbar$ if and only if
its characteristic polynomial has the form $x^3+\alpha$ for some $\alpha \in \K$, i.e., if and only
if $\tr(M)=c_2(M)=0$.
It follows that the $\overline{1}$-spec linear subspaces of $\Mat_3(\K)$ are the
totally isotropic subspaces
of $\frak{sl}_3(\K)$ for the symmetric bilinear form $(A,B) \mapsto \tr(AB)$.
Consider the non-degenerate symmetric bilinear form $b(A,B):=\tr(AB)$ on $\Mat_3(\K)$, and notice that
$\K I_3+\NT_3(\K)$ is a $4$-dimensional totally isotropic subspace of it. Since $4=[9/2]$, it follows that
the Witt index of $b$ is $4$. Since $\frak{sl}_3(\K)=\{I_3\}^\bot$ and $I_3$ is $b$-isotropic,
the hyperbolic inflation theorem yields that all the maximal totally isotropic subspaces of
$\frak{sl}_3(\K)$ have dimension $4$ (which gives us a new proof of Corollary \ref{majocor} in that case)
and Witt's extension theorem shows that these subspaces form a single orbit under the (natural) action of the orthogonal group of
$(c_2)_{|\frak{sl}_3(\K)}$. This however gives us little information on their orbits under \emph{conjugation}, which is the topic of our investigation.

\vskip 3mm
In our study, it will be quite helpful to think in terms of spaces of linear transformations
(rather than solely of matrices):

\begin{Def}
Let $E$ be a finite-dimensional vector space over $\K$. Denote by $\End(E)$ its vector space of linear endomorphisms. \\
A linear subspace $V$ of $\End(E)$ is called a \textbf{$\overline{1}$-spec} subspace when every element of $V$
has a sole eigenvalue in $\overline{\K}$. \\
Given a basis $\bfB$ of $E$ (with cardinality $n$), we denote by $\Mat_\bfB(V)$ the linear subspace
consisting of the matrices of $\Mat_n(\K)$ representing the elements of $V$ in the basis~$\bfB$.
\end{Def}

In the rest of the section, we set $E:=\K^3$.

\begin{Def}
Let $V$ be a $4$-dimensional $\overline{1}$-spec subspace of $\End(E)$. \\
A vector $x \in E \setminus \{0\}$ is said to be \textbf{good for} $V$ when
no element of $V$ has $\Vect(x)$ as its range.
\end{Def}

\noindent Proposition \ref{linelemma} thus implies that at least one non-zero vector of $E$ is good for $V$.

\subsection{Finishing the reduction of an arbitrary $4$-dimensional $\overline{1}$-spec subspace}\label{finishred}

\begin{Not}
For $\delta \in \K$, we set
$$\calF_\delta:=\Vect\Biggl(I_3,
\begin{bmatrix}
0 & 1 & 0 \\
0 & 0 & 0 \\
1 & 0 & 0
\end{bmatrix},
\begin{bmatrix}
0 & 0 & 0 \\
0 & 0 & 1 \\
\delta & 0 & 0
\end{bmatrix},
\begin{bmatrix}
1 & 0 & 1 \\
-1 & 0 & 0 \\
-1 & -\delta & -1
\end{bmatrix}\Biggr)$$
and
$$\calG_\delta:=\Vect\Biggl(I_3,
\begin{bmatrix}
0 & 1 & 0 \\
0 & 0 & 0 \\
1 & 0 & 0
\end{bmatrix},
\begin{bmatrix}
0 & 0 & 0 \\
0 & 0 & 1 \\
\delta & 0 & 0
\end{bmatrix},
\begin{bmatrix}
0 & 0 & 1 \\
-1 & 0 & 0 \\
0 & -\delta & 0
\end{bmatrix}\Biggr).$$
\end{Not}

Let $V$ be a $4$-dimensional $\overline{1}$-spec subspace of $\End(E)$.
We may find a vector $e_3 \in E \setminus \{0\}$ which is good for $V$. \\
In Sections \ref{setup} and \ref{specialmat}, we have shown that we may find two vectors $e_1$ and $e_2$ in $E$
such that $\bfB=(e_1,e_2,e_3)$ is a basis of $E$ and
$\Mat_\bfB(V)=\Vect(I_3,A_1,B_1,J)$, where
$$A_1=\begin{bmatrix}
0 & 1 & 0 \\
0 & 0 & 0 \\
a & \lambda & 0
\end{bmatrix} \quad \text{and} \quad
B_1=\begin{bmatrix}
0 & 0 & 0 \\
\mu & 0 & 1 \\
b & 0 & 0
\end{bmatrix}$$
and $J$ has the form $\begin{bmatrix}
0 & 0 & 1 \\
? & ? & 0 \\
? & ? & ?
\end{bmatrix}$.
Since $\tr(A_1B_1)=0$, we find $\mu=-\lambda$. \\
As in Paragraph \ref{lastconjugate}, we may then modify $e_1$
so as to have $\lambda=\mu=0$ in the new basis (which we still denote by $(e_1,e_2,e_3)$).

Replacing $J$ with $J+t.I_3$ for a well chosen $t \in \K$, we find in $V$ a matrix of the form
$$J'=\begin{bmatrix}
t & 0 & 1 \\
? & 0 & 0 \\
? & ? & ?
\end{bmatrix}.$$
Then $\tr(J')=0$, $\tr(A_1J')=0$, $\tr(B_1J')=0$ and $c_2(J')=0$ yield:
$$J'=\begin{bmatrix}
t & 0 & 1 \\
-a & 0 & 0 \\
-t^2 & -b & -t
\end{bmatrix}.$$
At this point, we need to distinguish between several cases:

\begin{itemize}
\item Assume $a \neq 0$.
Choose an arbitrary $\gamma \in \K \setminus \{0\}$, and set
$\delta:=\dfrac{ab}{\gamma^3}$ and $\bfB':=\bigl(\frac{1}{\gamma}\,e_1,\frac{a}{\gamma^2}\,e_2,e_3\bigr)$.
Then $\Mat_{\bfB'}(V)$ is
spanned by
$I_3$,
$\begin{bmatrix}
0 & 1 & 0 \\
0 & 0 & 0 \\
1 & 0 & 0
\end{bmatrix}$,
$\begin{bmatrix}
0 & 0 & 0 \\
0 & 0 & 1 \\
\delta & 0 & 0
\end{bmatrix}$,
$\begin{bmatrix}
(t/\gamma) & 0 & 1 \\
-1 & 0 & 0 \\
-(t/\gamma)^2 & -\delta & -(t/\gamma)
\end{bmatrix}$.
\begin{itemize}
\item If $t=0$, then we deduce that $\Mat_{\bfB'}(V)=\calG_\delta$.
\item If $t \neq 0$, then, by \emph{choosing} $\gamma=t$, we have $\delta=\frac{ab}{t^3}$ and
$\Mat_{\bfB'}(V)=\calF_\delta$.
\end{itemize}

\vskip 2mm

\item Assume $a=0$ and $b \neq 0$. Then $\Mat_{(e_3,e_1,e_2)}(V)$ is spanned by $I_3$,
$\begin{bmatrix}
0 & 0 & 0 \\
0 & 0 & 1 \\
0 & 0 & 0
\end{bmatrix}$, $\begin{bmatrix}
0 & 1 & 0 \\
0 & 0 & 0 \\
b^{-1} & 0 & 0
\end{bmatrix}$ and
$\begin{bmatrix}
-t & -t^2 & -b \\
1 & t & 0 \\
0 & 0 & 0
\end{bmatrix}$.
Multiplying the last matrix by $-\frac{1}{b}$ and summing it with $\dfrac{t}{b}\cdot I_3-\dfrac{t^2}{b}\cdot \begin{bmatrix}
0 & 1 & 0 \\
0 & 0 & 0 \\
b^{-1} & 0 & 0
\end{bmatrix}$, we find that $\Mat_{(e_3,e_1,e_2)}(V)$ is spanned by
$I_3$, $\begin{bmatrix}
0 & 1 & 0 \\
0 & 0 & 0 \\
b^{-1} & 0 & 0
\end{bmatrix}$, $\begin{bmatrix}
0 & 0 & 0 \\
0 & 0 & 1 \\
0 & 0 & 0
\end{bmatrix}$ and $\begin{bmatrix}
-(t/b) & 0 & 1 \\
-1/b & 0 & 0 \\
-t^2/b^2 & 0 & t/b
\end{bmatrix}$, and we deduce, as above, that $\Mat_{(e_1,e_2,e_3)}(V)$ is similar to $\calF_0$ (if $t\neq 0$) or to $\calG_0$ (if $t=0$).

\vskip 2mm

\item Assume that $a=b=0$ and $t \neq 0$.
Then
$$\Mat_{(e_1,e_2,e_3-\frac{e_1}{t})}(V)=\Vect\Biggl(I_3,\begin{bmatrix}
0 & 1 & 0 \\
0 & 0 & 0 \\
0 & 0 & 0
\end{bmatrix},\begin{bmatrix}
0 & 0 & 0 \\
0 & 0 & 1 \\
0 & 0 & 0
\end{bmatrix},\begin{bmatrix}
0 & 0 & 0 \\
0 & 0 & 0 \\
-1/t^2 & 0 & 0
\end{bmatrix}\Biggr)$$
and hence
$$\Mat_{(e_1,e_2,e_3)}(V) \sim \calI:=\K I_3 \oplus \Biggl\{\begin{bmatrix}
0 & x & 0 \\
0 & 0 & y \\
z & 0 & 0
\end{bmatrix} \mid (x,y,z)\in \K^3\Biggr\}.$$
This shows in particular that $\calF_0$ itself is similar to $\calI$:
indeed, for $P=\begin{bmatrix}
1 & 0 & 1 \\
1 & -1 & -1 \\
0 & 0 & -1
\end{bmatrix}$, we find that $P^{-1}\calF_0 P$
is spanned by $I_3$, $\begin{bmatrix}
-1 & -1 & 0 \\
0 & -1 & 1 \\
-1 & 0 & -1
\end{bmatrix}$, $\begin{bmatrix}
0 & 0 & 0 \\
0 & 0 & 1 \\
0 & 0 & 0
\end{bmatrix}$ and $\begin{bmatrix}
0 & 0 & 0 \\
0 & 0 & 1 \\
1 & 0 & 0
\end{bmatrix}$, and it is thus also spanned by
$I_3$, $\begin{bmatrix}
0 & 1 & 0 \\
0 & 0 & 0 \\
0 & 0 & 0
\end{bmatrix}$, $\begin{bmatrix}
0 & 0 & 0 \\
0 & 0 & 1 \\
0 & 0 & 0
\end{bmatrix}$ and $\begin{bmatrix}
0 & 0 & 0 \\
0 & 0 & 0 \\
1 & 0 & 0
\end{bmatrix}$.

\vskip 2mm

\item Finally, if $a=b=t=0$, then we readily have $V=\K I_3+\NT_3(\K)$.
\end{itemize}

\noindent We may summarize some of the above results as follows:

\begin{lemme}\label{lfin1}
Let $(e_1,e_2,e_3)$ be a basis of $E$. Assume that there exists
$(a,b,t)\in \K^3$ such that $a \neq 0$ and $\Mat_{(e_1,e_2,e_3)}(V)$ is spanned by the matrices
$I_3$,
$\begin{bmatrix}
0 & 1 & 0 \\
0 & 0 & 0 \\
a & 0 & 0
\end{bmatrix}$,
$\begin{bmatrix}
0 & 0 & 0 \\
0 & 0 & 1 \\
b & 0 & 0
\end{bmatrix}$ and $\begin{bmatrix}
t & 0 & 1 \\
-a & 0 & 0 \\
-t^2 & -b & -t
\end{bmatrix}$. \\
If $t \neq 0$, then there exists $(\lambda,\mu)\in (\K \setminus \{0\})^2$ such that
$\Mat_{(\lambda.e_1,\mu.e_2,e_3)}(V)=\calF_{ab/t^3}$. \\
If $t=0$, then, for every $\alpha \in \K \setminus \{0\}$, there exists $(\lambda,\mu)\in (\K \setminus \{0\})^2$
such that $\Mat_{(\lambda.e_1,\mu.e_2,e_3)}(V)=\calG_{\alpha^3 ab}$. \\
\end{lemme}

\begin{prop}
The subspaces $\calI$ and $\calF_0$ are similar.
\end{prop}

\subsection{Classification theorems}

We have just proven a good deal of the following theorem:

\begin{theo}\label{egaltheo3}
Assume $\car(\K)=3$, and let $V$ be a $4$-dimensional $\overline{1}$-spec linear subspace of $\Mat_3(\K)$.
Then exactly one of the following statements is true:
\begin{enumerate}[(i)]
\item $V$ is similar to $\K I_3+\NT_3(\K)$;
\item $V$ is similar to $\calF_\delta$ for some $\delta \in \K$;
\item $V$ is similar to $\calG_\delta$ for some $\delta \in \K$.
\end{enumerate}
\end{theo}

\noindent Conversely, a straightforward computation shows that
$\calF_\delta$ and $\calG_\delta$ are $4$-dimensional $\overline{1}$-spec subspaces of $\Mat_3(\K)$, for every $\delta \in \K$.

\begin{Not}
Denote by $\sigma : x \mapsto x^3$ the Frobenius automorphism of $\K$,
and set $j:=\sigma-\id$. Note that $j$ is an endomorphism of the group $(\K,+)$.
\end{Not}

\begin{Not}
We define the relation $\simt$ on $\K$ as follows:
$$x \simt \,y \; \underset{\text{def}}{\Leftrightarrow}\; \exists (a,b)\in (\K\setminus \{0\}) \times \K : \; x=a^3y+b^3.$$
This is obviously an equivalence relation on $\K$.
\end{Not}

\begin{theo}\label{egaltheo3prec}
Assume $\car(\K)=3$. Let $(\delta,\lambda)\in \K^2$. Then:
$$\calG_\delta \sim \calG_\lambda \; \Leftrightarrow \; \delta \simt \,\lambda
\qquad \text{and} \qquad
\calF_\delta \sim \calF_\lambda \; \Leftrightarrow \; \delta=\lambda \; \text{mod. $j(\K)$.}$$
\end{theo}

\begin{ex}
Assume $\K$ is an algebraically closed field of characteristic $3$.
Then $j$ and $\sigma$ are onto, and therefore there are exactly three similarity classes of $4$-dimensional $\overline{1}$-spec subspaces of $\Mat_3(\K)$:
those of $\K I_3+\NT_3(\K)$, $\calF_0 \sim \calI$ and $\calG_0$.
\end{ex}

\begin{ex}
Assume $\K$ is a finite field of characteristic $3$.
Then $\sigma$ is onto and hence $\simt$ has a sole equivalence class,
whereas $j(\K)$ is a subgroup of index $3$ of $(\K,+)$ (since the kernel of $j$ is the prime subfield of $\K$ and has therefore three elements).
Therefore, there are five similarity classes of $4$-dimensional $\overline{1}$-spec subspaces of $\Mat_3(\K)$.
\end{ex}

\vskip 3mm
\noindent In order to completely classify the $4$-dimensional $\overline{1}$-spec subspaces up to similarity,
what remains to prove is Theorem \ref{egaltheo3prec} and the uniqueness statement in Theorem \ref{egaltheo3}:
we achieve this is the next section.

\subsection{The uniqueness statements in Theorem \ref{egaltheo3}, and the proof of Theorem \ref{egaltheo3prec}}

\begin{lemme}
Let $\delta \in \K$. Then neither $\calF_\delta$ nor $\calG_\delta$
is similar to $\K I_3+\NT_3(\K)$.
\end{lemme}

\begin{proof}
In $\K I_3+\NT_3(\K)$, the set of singular matrices is the linear subspace $\NT_3(\K)$.
It thus suffices to find two singular matrices in $\calF_\delta$ (respectively, $\calG_\delta$)
some linear combination of which is non-singular. \\
The matrices $A=\begin{bmatrix}
0 & 1 & 0 \\
0 & 0 & 0 \\
1 & 0 & 0
\end{bmatrix}$ and $B=\begin{bmatrix}
0 & 0 & 0 \\
0 & 0 & 1 \\
\delta & 0 & 0
\end{bmatrix}$ are singular and belong both to $\calF_\delta$ and $\calG_\delta$.
Choosing $t \in \K \setminus \{0,-\delta\}$, we find that $t\,A+B$ is non-singular, which completes the proof.
\end{proof}

\noindent We turn to the study of possible similarities between spaces of type $\calF_\delta$ or $\calG_\delta$.

\begin{Def}
A basis $\bfB$ of $V$ is called \textbf{$V$-adapted} if
there exists $\delta \in \K$ such that $\Mat_\bfB(V)=\calF_\delta$ or $\Mat_\bfB(V)=\calG_\delta$:
notice then that $\delta$ is uniquely determined by $\bfB$ and that only one of the
conditions $\Mat_\bfB(V)=\calF_\delta$ and $\Mat_\bfB(V)=\calG_\delta$ holds.
\end{Def}

Let $V$ be a $4$-dimensional $\overline{1}$-spec subspace of $\End(E)$ which is
not similar to $\K I_3+\NT_3(\K)$.
Note that the third vector of a $V$-adapted basis is always good for $V$.

\vskip 2mm
Here is our strategy:
\begin{itemize}
\item given a vector $e_3 \in E$ which is good for $V$, determine the spaces of matrices associated to
the $V$-adapted bases with $e_3$ as the last vector;
\item then investigate what happens when the last vector of a $V$-adapted basis is modified by an ``elementary" operation.
\end{itemize}

The following lemma states that some sort of converse statement of Lemma \ref{lfin1} holds:

\begin{lemme}\label{lfin2}
Let $(e_1,e_2,e_3)$ be a $V$-adapted basis of $E$.
Assume that there exists $(\lambda,\mu,a,b,t)\in \K^5$ such that $\lambda \neq 0$, $\mu \neq 0$,
and $\Mat_{(\lambda.e_1,\mu.e_2,e_3)}(V)$ is spanned by $I_3$,
$\begin{bmatrix}
0 & 1 & 0 \\
0 & 0 & 0 \\
a & 0 & 0
\end{bmatrix}$,
$\begin{bmatrix}
0 & 0 & 0 \\
0 & 0 & 1 \\
b & 0 & 0
\end{bmatrix}$ and $\begin{bmatrix}
t & 0 & 1 \\
-a & 0 & 0 \\
-t^2 & -b & -t
\end{bmatrix}$. \\
It $t \neq 0$, then $\Mat_{(e_1,e_2,e_3)}(V)=\calF_{ab/t^3}$. \\
If $t=0$, then $\Mat_{(e_1,e_2,e_3)}(V)=\calG_{\alpha^3 ab}$ for some $\alpha \in \K \setminus \{0\}$.
\end{lemme}

\begin{proof}
Let $\delta \in \K$ be such that $\Mat_{(e_1,e_2,e_3)}(V)=\calF_\delta$ or $\Mat_{(e_1,e_2,e_3)}(V)=\calG_\delta$.
Note that $\Mat_{(e_1,e_2,e_3)}(V)$ is spanned by $I_3$,
$\begin{bmatrix}
0 & \frac{\lambda}{\mu} & 0 \\
0 & 0 & 0 \\
\frac{a}{\lambda} & 0 & 0
\end{bmatrix}$,
$\begin{bmatrix}
0 & 0 & 0 \\
0 & 0 & \mu \\
\frac{b}{\lambda} & 0 & 0
\end{bmatrix}$ and
$\begin{bmatrix}
t & 0 & \lambda \\
-\frac{a \mu}{\lambda} & 0 & 0 \\
-\frac{t^2}{\lambda} & -\frac{b}{\mu} & -t
\end{bmatrix}$, and hence by
$I_3$, $\begin{bmatrix}
0 & 1 & 0 \\
0 & 0 & 0 \\
\frac{a\mu}{\lambda^2} & 0 & 0
\end{bmatrix}$,
$\begin{bmatrix}
0 & 0 & 0 \\
0 & 0 & 1 \\
\frac{b}{\lambda \mu} & 0 & 0
\end{bmatrix}$ and
$\begin{bmatrix}
\frac{t}{\lambda} & 0 & 1 \\
-\frac{a \mu}{\lambda^2} & 0 & 0 \\
-\bigl(\frac{t}{\lambda}\bigr)^2 & -\frac{b}{\lambda \mu} & -\frac{t}{\lambda}
\end{bmatrix}$. \\
We deduce that $a\mu=\lambda^2$ and $\delta=\frac{b}{\lambda\mu}=\frac{ab}{\lambda^3}\cdot$ \\
This immediately yields the second result if $t=0$. If $t \neq 0$, then
$\frac{t}{\lambda}=1$ hence $\lambda=t$, which yields the first result.
\end{proof}

\begin{prop}\label{biglastprop}
Let $\bfB=(e_1,e_2,e_3)$ be a $V$-adapted basis.
Let $\delta \in \K$.
\begin{enumerate}[(a)]
\item Assume $\Mat_{\bfB}(V)=\calG_\delta$.
\begin{enumerate}[(i)]
\item For every $\lambda \in \K$ such that $\lambda \simt \delta$, there exists
a $V$-adapted basis $\bfB'$ with $e_3$ as the last vector such that $\Mat_{\bfB'}(V)=\calG_\lambda$.
\item For every $V$-adapted basis $\bfB'$ with $e_3$ as the last vector, there exists $\lambda \in \K$ such that
$\lambda \simt \delta$ and $\Mat_{\bfB'}(V)=\calG_\lambda$.
\end{enumerate}

\item Assume $\Mat_{\bfB}(V)=\calF_\delta$.
\begin{enumerate}[(i)]
\item For every $\lambda \in \K$ such that $\lambda = \delta$ mod. $j(\K)$, there exists
a $V$-adapted basis $\bfB'$ with $e_3$ as the last vector such that $\Mat_{\bfB'}(V)=\calF_\lambda$.
\item For every $V$-adapted basis $\bfB'$ with $e_3$ as the last vector, there exists $\lambda \in \K$ such that
$\lambda=\delta$ mod. $j(\K)$ and $\Mat_{\bfB'}(V)=\calF_\lambda$.
\end{enumerate}
\item For every $s \in \K$, there is a $V$-adapted basis of the form $(?,f_2,e_3)$ such that $\Vect(f_2,e_3)=\Vect(e_2+s e_1,e_3)$.
\end{enumerate}

\end{prop}

\begin{proof}
We start with a preliminary computation.
Let $(s,u,t)\in \K^3$ and set
$P:=\begin{bmatrix}
1 & s & 0 \\
0 & 1 & 0 \\
s & u & 1
\end{bmatrix}$, and $V_0$ the subspace of $\Mat_3(\K)$ spanned by $I_3$,
$A=\begin{bmatrix}
0 & 1 & 0 \\
0 & 0 & 0 \\
1 & 0 & 0
\end{bmatrix}$,
$\begin{bmatrix}
0 & 0 & 0 \\
0 & 0 & 1 \\
\delta & 0 & 0
\end{bmatrix}$ and
$\begin{bmatrix}
t & 0 & 1 \\
-1 & 0 & 0 \\
-t^2 & -\delta & -t
\end{bmatrix}$. \\
A straightforward computation shows that $P$ is non-singular and $P^{-1}V_0 P$ is spanned by
$$I_3, \; A, \; \begin{bmatrix}
-s^2 & -us & -s \\
s & u & 1 \\
\delta+s^3-us & ? & -u+s^2
\end{bmatrix} \; \text{and} \;
\begin{bmatrix}
t-s & ts+s^2+u & 1 \\
-1 & -s & 0 \\
ts+s^2+u-t^2 & ? & -s-t
\end{bmatrix}.$$
Adding well-chosen linear combinations of the first two matrices to the last two, we find that $P^{-1}V_0 P$ is spanned by
$$I_3, \; A, \; \begin{bmatrix}
0 & 0 & -s \\
s & u+s^2 & 1 \\
\delta+s^3 & ? & -u-s^2
\end{bmatrix} \; \text{and} \;
\begin{bmatrix}
0 & 0 & 1 \\
-1 & -t & 0 \\
-t^2 & ? & t
\end{bmatrix}.$$
Adding to the third one the product of the fourth one with $s$, we deduce that
$P^{-1}V_0 P$ is spanned by
$$I_3, \; A, \; \begin{bmatrix}
0 & 0 & 0 \\
0 & u+s^2-st & 1 \\
\delta+s^3-st^2 & -(u+s^2-st)^2 & -u-s^2+st
\end{bmatrix} \; \text{and} \;
\begin{bmatrix}
0 & 0 & 1 \\
-1 & -t & 0 \\
-t^2 & ? & t
\end{bmatrix}$$
(the $(3,2)$-th entry of the third matrix is easily obtained using the fact that this matrix is nilpotent).
Letting $(g_1,g_2,g_3)$ be an arbitrary basis of $\K^3$ and denoting by $H$ the linear subspace of $\End(\K^3)$ such that
$\Mat_{(g_1,g_2,g_3)}(H)=V_0$, we then find that
$u=st-s^2$ if and only if there exists $(a,b)\in (\K \setminus \{0\})^2$ such that
$\bigl(a.(g_1+s.g_3),b.(g_2+s.g_1+u.g_3),g_3\bigr)$ is $H$-adapted. \\
If $u=st-s^2$, then
$P^{-1}V_0 P$ is spanned by
$$I_3, \; A, \; \begin{bmatrix}
0 & 0 & 0 \\
0 & 0 & 1 \\
\delta+s^3-st^2 & 0 & 0
\end{bmatrix} \; \text{and} \;
\begin{bmatrix}
t & 0 & 1 \\
-1 & 0 & 0 \\
-t^2 & st^2-s^3-\delta & -t
\end{bmatrix},$$
(the $(3,2)$-th entry of the fourth matrix is obtained using the fact that the third and fourth matrices are mutually orthogonal for
$(M,N) \mapsto \tr(MN)$)
and we deduce that $P^{-1}V_0 P=\calG_{\delta+s^3}$ if $t=0$, and $P^{-1}V_0 P=\calF_{\delta+s^3-s}$ if $t=1$.

\begin{itemize}
\item Proof of (a)(i):
Let $s \in \K$. Then the above calculation shows that
$\Mat_{(e_1+s e_3,e_2+s e_1-s^2 e_3,e_3)}(V)=\calG_{\delta+s^3}$.
Using Lemma \ref{lfin1}, we deduce that, for every $(s,z)\in \K \times (\K \setminus \{0\})$,
there exists a $V$-adapted basis $\bfB'$ with $e_3$ as third vector such that
$\Mat_{\bfB'}(V)=\calG_{z^3 \delta+(zs)^3}$. This obviously yields point (i).

\item Proof of (a)(ii): Let $(f_1,f_2)\in E^2$ be such that $(f_1,f_2,e_3)$ is a $V$-adapted basis. \\
Then the set of elements of $V$ which vanish on $\Vect(e_3)$ must be a $1$-dimensional subspace.
Given such a non-zero element $u$, we then have $\im(u)=\Vect(f_1,e_3)$. \\
It follows that $\Vect(e_1,e_3)=\Vect(f_1,e_3)$, which proves that $f_1=a(e_1+s e_3)$
for some $a \in \K \setminus \{0\}$ and some $s \in \K$. \\
Moreover, $u$ may be chosen so as to have $u(e_2)=e_1$, and there exists
$b \in \K \setminus \{0\}$ such that $u(f_2)=b\, f_1$ and $u(f_1)=b\, e_3$.
Therefore $u(f_2)=ab(e_1+s e_3)=ab \,u(e_2+s e_1)$, which yields a scalar
$\mu \in \K$ such that $f_2=ab\,(e_2+s e_1+\mu e_3)$. \\
Set now $e'_1:=e_1+s e_3$ and $e'_2:=e_2+s e_1+\mu e_3$. \\
Then the preliminary calculations show that $\mu=-s^2$ and
$\Mat_{(e_1+s e_3,e_2+s e_1+\mu e_3,e_3)}(V)=\calG_{\delta+s^3}$.
Using Lemma \ref{lfin2}, we deduce that $\Mat_{(f_1,f_2,e_3)}(V)=\calG_\lambda$ for some
$\lambda \in \K$ satisfying $\lambda \simt \delta$.
\end{itemize}
Points (b)(i) and (b)(ii) are deduced from the preliminary computation and Lemmas \ref{lfin1} and \ref{lfin2} in the same fashion
by taking $t=1$. \\
Point (c) follows directly from the proofs of points (a)(i) and (b)(i).
\end{proof}

We now examine the effect of a simple change of the last vector of a $V$-adapted basis.

\begin{lemme}\label{lfin3}
Let $(e_1,e_2,e_3)$ be a $V$-adapted basis of $E$.
Let $\delta \in \K$ such that $\Mat_{(e_1,e_2,e_3)}(V)=\calG_\delta$ (respectively, $\Mat_{(e_1,e_2,e_3)}(V)=\calF_\delta$).
\begin{enumerate}[(a)]
\item Assume that $\delta \underset{3}{\not\sim} 0$ (respectively, $\delta \neq 0$ mod. $j(\K)$). \\
Then, for every $s \in \K$, the vector $e_3+s e_2$ is the third one of a
$V$-adapted basis $\bfB$ such that $\Mat_{\bfB}(V)=\calG_\delta$ (respectively, $\Mat_{\bfB}(V)=\calF_\delta$).
\item If $\delta=0$, then
$e_3+e_2$ is the third vector of a $V$-adapted basis $\bfB$ such that $\Mat_{\bfB}(V)=\calG_\delta$ (respectively, $\Mat_{\bfB}(V)=\calF_\delta$).
\end{enumerate}
\end{lemme}

\begin{proof}
Set $t:=0$ if $\Mat_{(e_1,e_2,e_3)}(V)=\calG_\delta$, and $t:=1$ if $\Mat_{(e_1,e_2,e_3)}(V)=\calF_\delta$. \\
Let $s \in \K$.
Set $e'_1:=e_1-s(t+\delta s)e_2+\delta s e_3$, and $e'_3:=e_3+s e_2$. \\
Then a straightforward computation shows that $\Mat_{(e'_1,e_2,e'_3)}(V)$ is spanned by the matrices
$$I_3, \; \begin{bmatrix}
-s(t+\delta s) & 1 & s \\
-s-s^2t^2+s^4\delta^2  & s(t-\delta s) & s^2(t-\delta s) \\
1+\delta s^2(t-\delta s)-\delta^2 s^3 & -\delta s & -\delta s^2
\end{bmatrix}, \;
\begin{bmatrix}
0 & 0 & 0 \\
0 & 0 & 1 \\
\delta & 0 & 0
\end{bmatrix} \; \text{and} \;
\begin{bmatrix}
\delta s+t & 0 & 1 \\
-st^2+\delta^2 s^3-1 & \delta s & -st \\
-\delta s t-t^2 & -\delta & \delta s-t
\end{bmatrix}.$$
Adding well-chosen linear combinations of $I_3$ and the third matrix to the second and the fourth ones, we deduce that
$\Mat_{(e'_1,e_2,e'_3)}(V)$ is spanned by the matrices
$$I_3, \; A_1=\begin{bmatrix}
0 & 1 & s \\
-s-s^2t^2+s^4\delta^2  & -st & 0 \\
1-\delta^2 s^3 & -\delta s & st
\end{bmatrix}, \;
A_2=\begin{bmatrix}
0 & 0 & 0 \\
0 & 0 & 1 \\
\delta & 0 & 0
\end{bmatrix} \; \text{and} \;
A_3=\begin{bmatrix}
0 & 0 & 1 \\
-s t^2+\delta^2 s^3-1 & -t & 0 \\
-t^2 & -\delta & t
\end{bmatrix}.$$
By $A_1 \leftarrow A_1-sA_3$, we deduce that
$\Mat_{(e'_1,e_2,e'_3)}(V)$ is spanned by the matrices
$$I_3, \; \begin{bmatrix}
0 & 1 & 0 \\
0  & 0 & 0 \\
1-\delta^2 s^3+s t^2 & 0 & 0
\end{bmatrix}, \;
\begin{bmatrix}
0 & 0 & 0 \\
0 & 0 & 1 \\
\delta & 0 & 0
\end{bmatrix} \; \text{and} \;
\begin{bmatrix}
0 & 0 & 1 \\
-1+\delta^2 s^3-s t^2 & -t & 0 \\
-t^2 & -\delta & t
\end{bmatrix}.$$
Set $\gamma:=\delta(1-\delta^2 s^3+s t^2)$ and notice that
$\gamma=\delta+(-\delta s)^3$ if $t=0$, and
$\gamma=\delta+j(-\delta s)$ if $t=1$. \\
In any case, the assumptions of point (a) show that $\gamma \neq 0$ and therefore $1-\delta^2 s^3+s t^2 \neq 0$:
Lemmas \ref{lfin1} and \ref{lfin2} thus yield a $V$-adapted basis $\bfB$ with $e'_3$ as third vector such that
$\Mat_\bfB(V)=\calG_{\delta}$ if $t=0$, and $\Mat_\bfB(V)=\calF_{\delta}$ if $t=1$. \\
Assume now that $\delta=0$ and $s=1$. Then $1-\delta^2 s^3+s t^2 \neq 0$, hence Lemma \ref{lfin1}
yields point (b).
\end{proof}

\begin{lemme}
Let $(e_1,e_2,e_3)$ be a $V$-adapted basis of $E$.
Let $\delta \in \K$ such that $\Mat_{(e_1,e_2,e_3)}(V)=\calG_\delta$ (respectively, $\Mat_{(e_1,e_2,e_3)}(V)=\calF_\delta$),
and assume that $\delta \neq 0$. Then $e_2$ is the third vector of a $V$-adapted basis $\bfB$
such that $\Mat_\bfB(V)=\calG_\delta$ (respectively, $\Mat_\bfB(V)=\calF_\delta$).
\end{lemme}

\begin{proof}
Set $\bfB:=(e_3,e_1,e_2)$, and $t:=0$ if $\Mat_{(e_1,e_2,e_3)}(V)=\calG_\delta$,
and $t:=1$ if $\Mat_{(e_1,e_2,e_3)}(V)=\calF_\delta$. \\
Then a straightforward computation shows that
$\Mat_{(e_3,e_1,e_2)}(V)$ is spanned by
$$I_3, \; A_1=\begin{bmatrix}
0 & 1 & 0 \\
0 & 0 & 1 \\
0 & 0 & 0
\end{bmatrix}, \; A_2=\begin{bmatrix}
0 & \delta & 0 \\
0 & 0 & 0 \\
1 & 0 & 0
\end{bmatrix} \; \text{and} \; A_3=\begin{bmatrix}
-t & -t^2 & -\delta \\
1 & t & 0 \\
0 & -1 & 0
\end{bmatrix}.$$
Setting $A'_2:=\frac{1}{\delta}\,A_2$ and then $A'_1:=A_1-A'_2$, we find by
$\Mat_{(e_3,e_1,e_2)}(V)$ is spanned by
$$I_3, \; A'_1=\begin{bmatrix}
0 & 0 & 0 \\
0 & 0 & 1 \\
-\delta^{-1} & 0 & 0
\end{bmatrix}, \; A'_2=\begin{bmatrix}
0 & 1 & 0 \\
0 & 0 & 0 \\
\delta^{-1} & 0 & 0
\end{bmatrix} \; \text{and} \; A_3=\begin{bmatrix}
-t & -t^2 & -\delta \\
1 & t & 0 \\
0 & -1 & 0
\end{bmatrix}.$$
Setting finally $A'_3:=-\delta^{-1}(A_3-t\,I_3+t^2\,A'_2)$, we find that
$\Mat_{(e_3,e_1,e_2)}(V)$ is spanned by
$I_3$, $A'_1$, $A'_2$ and the matrix
$$A'_3=\begin{bmatrix}
-\frac{t}{\delta} & 0 & 1 \\
-\delta^{-1} & 0 & 0 \\
-\bigl(\frac{t}{\delta}\bigr)^2 & \delta^{-1} & \frac{t}{\delta}
\end{bmatrix}.$$
Assume $t=0$. Notice that $\delta=(-\delta)^3 \times (-\delta^{-1}) \times \delta^{-1}$, hence
Lemma \ref{lfin1} shows that $e_2$ is the third vector of a basis $\bfB$ such that
$\Mat_\bfB(V)=\calG_\delta$. \\
Assume $t=1$. Notice that $\frac{1}{(-t/\delta)^3}\,(-\delta^{-1}) \times \delta^{-1}=\delta$,
hence Lemma \ref{lfin1} shows that $e_2$ is the third vector of a basis $\bfB$ such that
$\Mat_\bfB(V)=\calF_\delta$.
\end{proof}

We are now ready to conclude. Let us first draw a corollary from the two previous lemmas:

\begin{cor}\label{lastcor}
Let $(e_1,e_2,e_3)$ be a $V$-adapted basis of $E$.
Let $\delta \in \K$ such that $\Mat_{(e_1,e_2,e_3)}(V)=\calG_\delta$ (respectively, $\Mat_{(e_1,e_2,e_3)}(V)=\calF_\delta$),
and assume that $\delta \neq 0$. Then, for any non-zero vector $z$ of $\Vect(e_2,e_3)$,
there is a $V$-adapted basis $\bfB$ with $z$ as the last vector
such that $\Mat_\bfB(V)=\calG_\delta$ (respectively, $\Mat_\bfB(V)=\calF_\delta$).
\end{cor}

This uses of course the trivial fact that given an arbitrary basis $(f_1,f_2,f_3)$ of $E$,
one has $\Mat_{(\lambda f_1,\lambda f_2, \lambda f_3)}(V)=\Mat_{(f_1,f_2,f_3)}(V)$ for any $\lambda \in \K \setminus \{0\}$.

\paragraph{}
We now prove Theorem \ref{egaltheo3prec} and the fact that an $\calF_\delta$-space is never similar to a $\calG_\lambda$-space.
First, the implications ``$\Leftarrow$" in Theorem \ref{egaltheo3prec} both follow directly from points (a)(i) and (b)(i) of Proposition \ref{biglastprop}. \\
For the converse implications,
let $V$ be a $4$-dimensional $\overline{1}$-spec linear subspace of $\End(E)$, where $E=\K^3$,
and assume that there exists a $V$-adapted basis $(e_1,e_2,e_3)$
with $\Mat_{(e_1,e_2,e_3)}(V)=\calG_\delta$ (respectively, $\Mat_{(e_1,e_2,e_3)}(V)=\calF_\delta$).
Of course, we may assume that $\delta \underset{3}{\not\sim} 0$ (respectively, $\delta \not\in j(\K)$). \\
Using point (c) of Proposition \ref{biglastprop} together with Corollary \ref{lastcor}, we find that every vector of $E \setminus \Vect(e_1,e_3)$
is the third vector of a $V$-adapted basis $\bfB$ for which $\Mat_\bfB(V)=\calG_\delta$ (respectively, $\Mat_\bfB(V)=\calF_\delta$).
Choose such a basis $(e'_1,e'_2,e'_3)$. Then
every vector of $E \setminus \Vect(e'_1,e'_3)$ is the third vector of a $V$-adapted basis $\bfB$ for which $\Mat_\bfB(V)=\calG_\delta$
(respectively, $\Mat_\bfB(V)=\calF_\delta$). However $\Vect(e'_1,e'_3) \cap \Vect(e_1,e_3)=\Vect(y)$ for some $y \in E \setminus \{0\}$
(because $e'_3 \not\in \Vect(e_1,e_3)$ and $\dim E=3$). \\
Therefore:
\begin{center}
(P) : \quad Every vector of $E \setminus \Vect(y)$ is the third vector of a $V$-adapted basis $\bfB$ for which $\Mat_\bfB(V)=\calG_\delta$
(respectively, $\Mat_\bfB(V)=\calF_\delta$).
\end{center}
By Lemma \ref{lfin2}, for every $V$-adapted basis $\bfB$ with a third vector which is linearly independent from $y$,
there exists $\lambda \in \K$ such that $\Mat_\bfB(V)=\calG_\lambda$ and $\lambda \simt \delta$
(respectively, $\Mat_\bfB(V)=\calF_\lambda$ and $\lambda=\delta$ mod. $j(\K)$). \\
Assume that $y$ is the third vector of a $V$-adapted basis $(g_1,g_2,y)$.
If $\Mat_{(g_1,g_2,y)}(V)$ equals $\calF_\lambda$ for some $\lambda \in j(\K)$ or
$\calG_\lambda$ for some $\lambda \in \sigma(\K)$, then Lemma \ref{lfin1} shows that
we lose no generality in assuming that $\lambda=0$,
and point (b) in Lemma \ref{lfin3} yields a vector $z\in E \setminus \Vect(y)$
which is the last vector of a $V$-adapted basis $\bfB'$ for which $\Mat_{\bfB'}(V)=\calG_0$
or $\Mat_{\bfB'}(V)=\calF_0$: this would contradict property (P). \\
Hence $\Mat_{(g_1,g_2,y)}(V)$ equals $\calF_\lambda$ for some $\lambda \in \K \setminus j(\K)$ or
$\calG_\lambda$ for some $\lambda \in \K \setminus \sigma(\K)$.
Applying point (a) of Proposition \ref{biglastprop}, we find a $V$-adapted basis $\bfB'$ with a third vector outside of
$\Vect(y)$ such that $\Mat_{\bfB'}(V)=\calF_\lambda$ or $\Mat_{\bfB'}(V)=\calG_\lambda$:
with the previous results, we deduce that $\Mat_{\bfB'}(V)=\calG_\lambda$ and
$\lambda \simt \delta$ (respectively, $\Mat_{\bfB'}(V)=\calF_\lambda$ and
$\lambda=\delta$ mod. $j(\K)$). \\
We have proven that, for every $V$-adapted basis, the subspace of matrices represented by $V$ in this basis
is $\calG_\lambda$ for some $\lambda \simt \delta$ (respectively, $\calF_\lambda$ for some $\lambda = \delta$ mod. $j(\K)$). \\
This finishes the proof of Theorem \ref{egaltheo3prec}.

We have also proved that, for any $\delta \in \K \setminus \sigma(\K)$,
the space $\calG_\delta$ is never similar to an $\calF_\lambda$-space, and, for any
$\delta \in \K \setminus j(\K)$, the space $\calF_\delta$ is never similar to a $\calG_\lambda$-space. \\
In order to complete the proof of Theorem \ref{egaltheo3}, it only remains to establish the following result:

\begin{lemme}
The spaces $\calF_0$ and $\calG_0$ are not similar.
\end{lemme}

\begin{proof}
We have already shown that $\calF_0$ is similar to the space
$$\K I_3+\Biggl\{\begin{bmatrix}
0 & x & 0 \\
0 & 0 & y \\
z & 0 & 0
\end{bmatrix}\mid (x,y,z)\in \K^3\Biggr\},$$
which obviously contains two linearly independent rank $1$ matrices. \\
We prove that $\calG_0$ does not contain such matrices. \\
Let $M \in \calG_0$ with rank $1$. We write $M=\begin{bmatrix}
\lambda & x & z \\
-z & \lambda & y \\
x &  0  & \lambda
\end{bmatrix}$.
Considering the lower right $2\times 2$ sub-matrix, we deduce from $\rk M \leq 1$ that $\lambda=0$. \\
Considering the sub-matrix obtained by deleting the second row and the third column, we find $x=0$.
It easily follows that $z=0$, therefore $M=y\,\begin{bmatrix}
0 & 0 & 0 \\
0 & 0 & 1 \\
0 &  0 & 0
\end{bmatrix}$. Two rank $1$ matrices of $\calG_0$ must therefore be linearly dependent, which
shows that $\calG_0$ is not similar to $\calF_0$.
\end{proof}

This finishes the proof of Theorem \ref{egaltheo3}. We have therefore completely classified
the $4$-dimensional $\overline{1}$-spec subspaces of $\Mat_3(\K)$, up to similarity.

\end{document}